\documentclass[lineno]{biometrika2}

\pdfoutput=1

\usepackage{amsmath}
\usepackage{array}
\usepackage[english]{babel}
\usepackage{appendix}
\usepackage{graphicx}
\usepackage{float}
\usepackage{bbm}
\usepackage{amsfonts}
\usepackage{mathbbol}
\usepackage{latexsym}
\usepackage{url}
\usepackage[T1]{fontenc}
\usepackage[utf8]{inputenc}
\usepackage{gensymb}

\usepackage{times}
\usepackage{bm}
\usepackage{natbib}

\usepackage[plain,noend]{algorithm2e}

\makeatletter
\renewcommand{\algocf@captiontext}[2]{#1\algocf@typo. \AlCapFnt{}#2} 
\def\@algocf@capt@plain{top}
\renewcommand{\algocf@makecaption}[2]{%
  \addtolength{\hsize}{\algomargin}%
  \sbox\@tempboxa{\algocf@captiontext{#1}{#2}}%
  \ifdim\wd\@tempboxa >\hsize
    \hskip .5\algomargin%
    \parbox[t]{\hsize}{\algocf@captiontext{#1}{#2}}
  \else%
    \global\@minipagefalse%
    \hbox to\hsize{\box\@tempboxa}
  \fi%
  \addtolength{\hsize}{-\algomargin}%
}
\makeatother

\begin{document}

\jname{}
\jyear{}
\jvol{}
\jnum{}
\copyrightinfo{}



\markboth{J. M{\o}ller, F. Safavimanesh, and J. G. Rasmussen}{The cylindrical $K$-function and Poisson line cluster point processes}

\title{The cylindrical $K$-function and Poisson line cluster point processes}

\author{Jesper M{\o}ller}
\affil{Department of Mathematical Sciences, Aalborg University, 9220 Aalborg, Denmark \email{jm@math.aau.dk}}
\author{Farzaneh Safavimanesh}
\affil{Department of Statistics, Faculty of Mathematical Sciences, Shahid Beheshti University, 19834 Tehran, Iran \email{f\_safavimanesh@sbu.ac.ir}}
\author{\and Jakob Gulddahl Rasmussen}
\affil{Department of Mathematical Sciences, Aalborg University, 9220 Aalborg, Denmark \email{jgr@math.aau.dk}}

\maketitle

\thispagestyle{empty}

\begin{abstract}
  The analysis of point patterns with linear structures is of interest
  in many applications.  To detect anisotropy in such cases, in
  particular in case of a columnar structure, we introduce a
  functional summary statistic, the cylindrical $K$-function, which is
  a directional $K$-function whose structuring element is a
  cylinder. Further we introduce a class of anisotropic Cox point
  processes, called Poisson line cluster point processes. The points
  of such a process are random displacements of Poisson point
  processes defined on the lines of a Poisson line process. Parameter
  estimation based on moment methods or Bayesian inference for this
  model is discussed when the underlying Poisson line process is
  latent.  To illustrate the methodologies, we analyze two- and
  three-dimensional point pattern data sets. The three-dimensional
  data set is of particular interest as it relates to the minicolumn
  hypothesis in neuroscience, claiming that pyramidal and other brain
  cells have a columnar arrangement perpendicular to the surface of
  the brain.
\end{abstract}


\section{Introduction}\label{intro}

Frequently in the spatial point process literature, isotropy, i.e.,
distributional invariance under rotations about a fixed location in
space, is assumed, though it is often unrealistic. Anisotropy of
spatial point processes has usually been studied by summarizing the
information of observed pairs of points, including the use of
directional $K$-functions or related densities
\citep{OhSt-81,StoyanBenes-91,Stoyan-91,StoyStoy-94,guan-etal:06,Illianetal-08,
  Redenbachetal-09}, spectral and wavelet methods
\citep*{MuggRen-96,Rosenberg-04,Nicolisetal-10}, and geometric
anisotropic pair correlation functions \citep{MoellerHokan-14}. The
applications considered in these references except
\citet{Redenbachetal-09} and \citet{Illianetal-08} are for
two-dimensional but not three-dimensional point patterns.

There are point patterns where points lie approximately along straight
lines, cf.\ the examples of applications and references in
\citet{moeller:waagepetersen:16}, called point patterns with linear
structures.  This paper focuses on detecting and modelling such point
patterns observed within a bounded subset of $\mathbb R^d$, $d\ge 2$,
where the cases $d=2$ and $d=3$ are of main interest. In particular we
study columnar structures. 

\begin{figure}[t]
\centering
\includegraphics[width=0.75\textwidth]{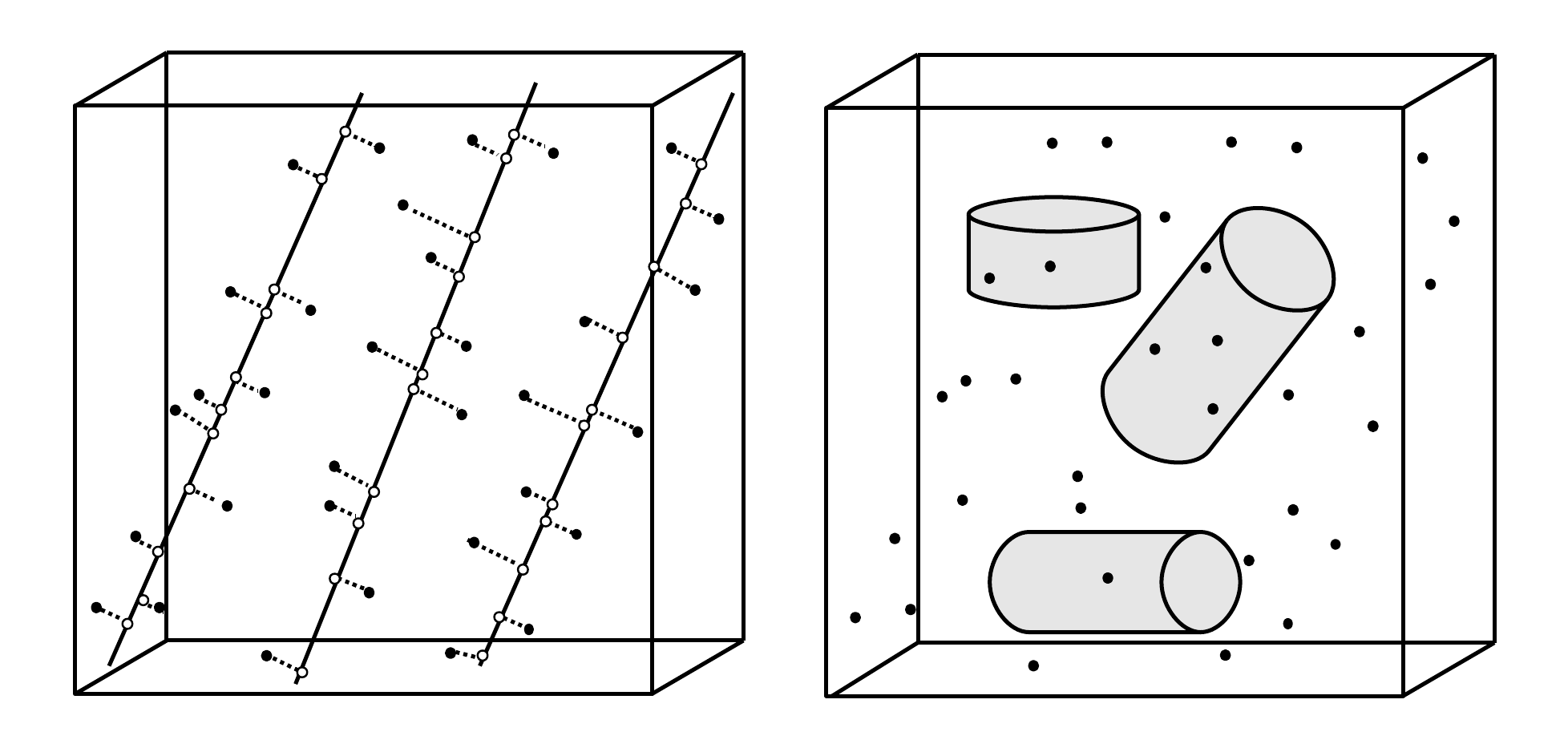}
\caption{Poisson line cluster process: Left panel: A simulated
  realization of a Poisson line cluster point process within a
  three-dimensional box.  The realizations of the Poisson line process
  (solid lines) and the Poisson point processes on the lines (circles)
  are also shown.  The dotted lines indicate how the points on the
  lines have been displaced to new positions (filled circles) and they
  specify the clusters.  Right panel: The same simulated realization
  of a Poisson line cluster point process and different choices of
  cylinders centered at different points of the process.}
\label{fig:cylinder}
\end{figure}

Section~\ref{s:cylK} introduces the cylindrical $K$-function, a
directional $K$-function whose structuring element is a cylinder which
is suitable for detecting anisotropy caused by columnar or other
linear structures in spatial point patterns. This is an adapted
version of the space-time $K$-function
\citep{Diggleetal-95,GabrielDiggel-09}.  Section~\ref{s:PLCPP}
concerns a new class of point processes, Poisson line cluster point
processes, whose points cluster around a Poisson line process. The
left panel in Figure~\ref{fig:cylinder} illustrates how such a process
is constructed: lines are generated from an anisotropic Poisson line
process, independent stationary Poisson point processes are generated
on the lines, and their points are randomly displaced, resulting in
the Poisson line cluster point process. We consider the Poisson lines
and the points on the lines as latent, so the clusters of the Poisson
line cluster point process are also hidden. Section~\ref{s:PLCPP} also
discusses a moment based approach and a simulation-based Bayesian
approach for inference, where in the latter case we estimate both the
parameters of the model and the missing lines.

\begin{figure}[t]
\centering
\includegraphics[width=7cm]{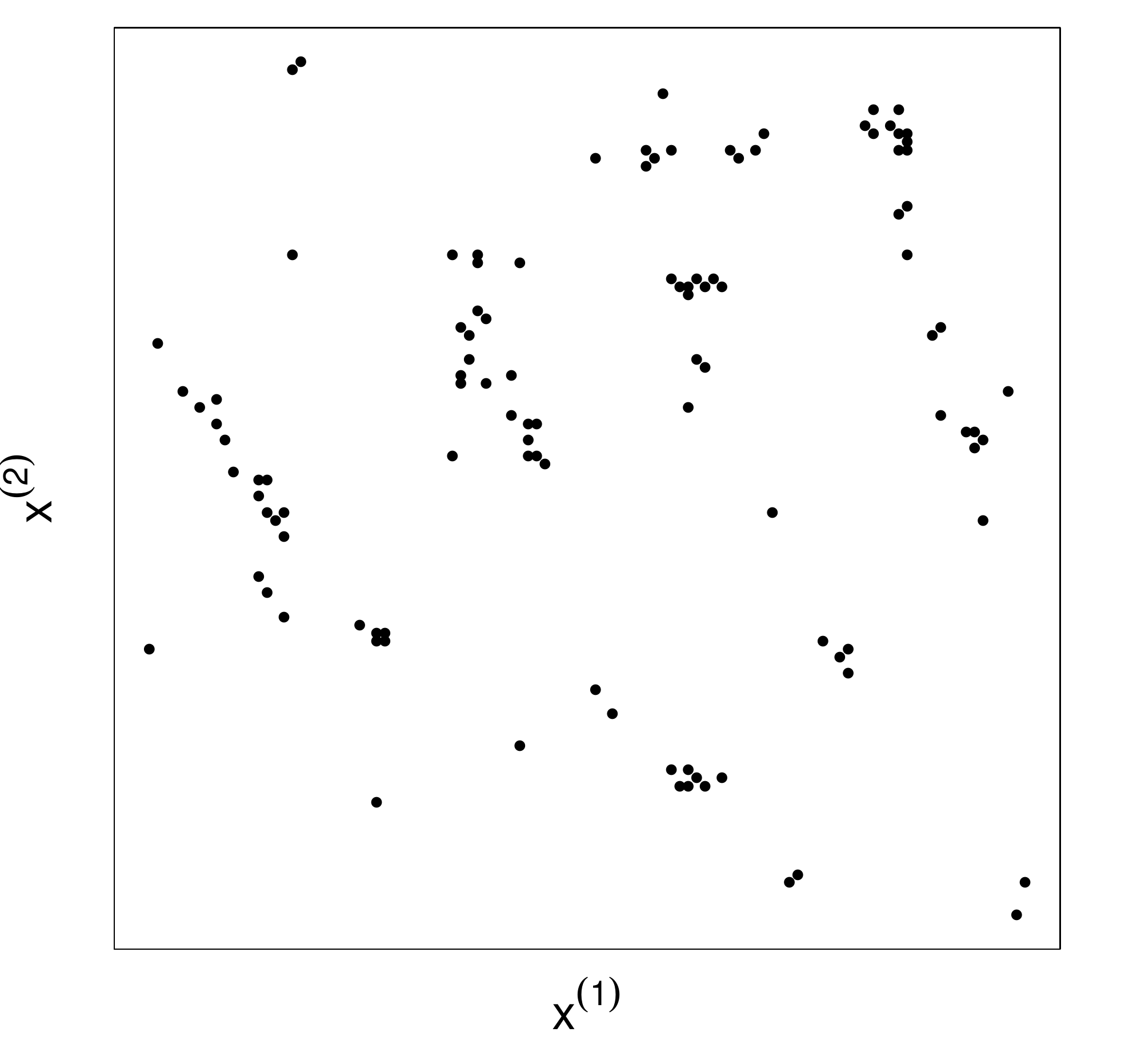}
\includegraphics[width=7cm]{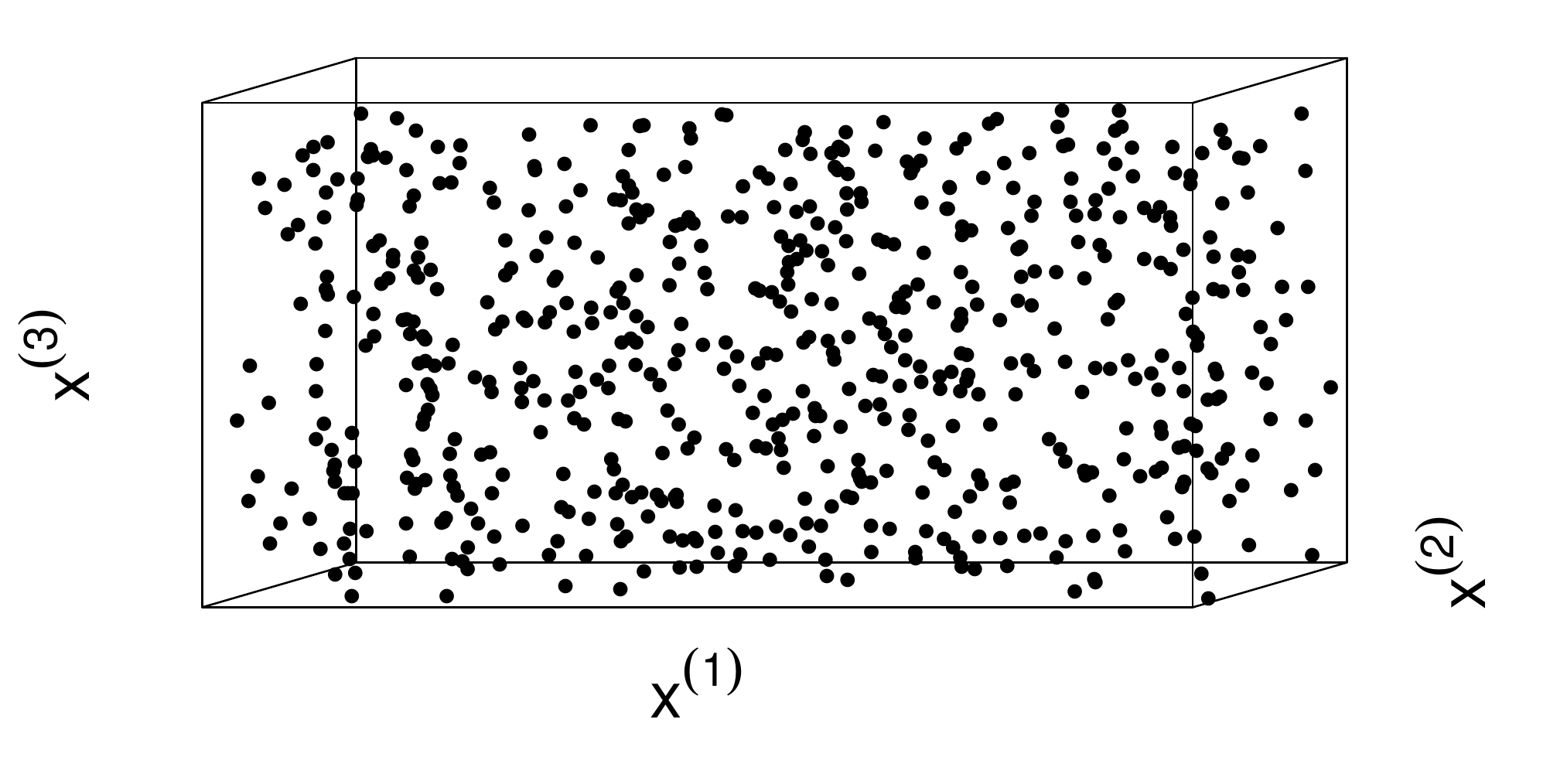}
\caption{Data sets: Left panel: Locations of 110 chapels in Wales,
  United Kingdom, observed in a square window (normalized to a unit
  square). Right panel: Nucleolus of 623 pyramidal cells in an
  observation window of size $508\times 138\times 320\,\mu{m}^3$.}
\label{fig:Data}
\end{figure}

Sections~\ref{s:cylK}-\ref{s:PLCPP} apply our methodology to the data
sets in Figure~\ref{fig:Data}.  The left panel shows a two-dimensional
point pattern data set recorded by \citet{MuggRen-96}, namely the
locations of 110 chapels in the Welsh Valleys, United Kingdom, where
the clear linear orientation is caused by four more or less parallel
valleys. For a three-dimensional point pattern data set it is often
difficult to detect anisotropy by eye, but the cylindrical
$K$-function will be useful, as illustrated later in connection to the
right panel, which shows the locations of 623 pyramidal cells from the
Brodmann area 4 of the grey matter of the human brain collected by
the neuroscientists at the Center for Stochastic Geometry and
Bioimaging, Denmark.  According to the minicolumn hypothesis
\citep{Mountcastle-57}, brain cells, mainly pyramidal cells, should
have a columnar arrangement perpendicular to the pial surface of the
brain, i.e., a columnar arrangement parallel to the $x^{(3)}$-axis
indicated in Figure~\ref{fig:Data}, and this should be highly
pronounced in Brodmann area 4. However, this hypothesis has been much
debated, see \citet{RSDRMN-15} and the references therein.

We use the chapel data set mainly for illustrative purposes and for
comparison with previous work. Investigations of the minicolumn
hypothesis have so far only been done in two dimensions except for the
three-dimensional analysis in \citet{RSDRMN-15}.  The present paper
details the methodology and provides a more thorough analysis of the
pyramidal cell data set, shedding further light on the validity of the
minicolumn hypothesis.

Throughout Sections~\ref{s:cylK}--\ref{s:PLCPP} we assume
stationarity. Section~\ref{s:extensions} discusses the choice of a
cylinder as the structuring element of the $K$-function
and extensions of this function and the Poisson line cluster point
process in a non-stationary setting.

\section{The cylindrical $K$-function}\label{s:cylK}

\subsection{Setting}\label{s:setting}

Throughout this paper we make the following assumptions and use the
following notation. 

Unless otherwise stated, we consider a stationary point process $X$
defined on $\mathbb R^d$, with finite and positive intensity $\rho$,
and where we view $X$ as a locally finite random subset of $\mathbb
R^d$.  Here stationarity means that the distribution of $X$ is
invariant under translations in $\mathbb R^d$, and $\rho|B|$ is the
mean number of points from $X$ falling in any Borel set
$B\subseteq\mathbb R^d$ of volume $|B|$.  We assume that $X$ has a
pair correlation function $g(x)$ defined for all $x\in\mathbb R^d$.
Intuitively, if $x_1,x_2\in\mathbb R^d$ are distinct locations and
$B_1,B_2$ are infinitesimally small sets of volumes $\mathrm
dx_1,\,\mathrm dx_2$ and containing $x_1,x_2$, respectively, then
$\rho^2g(x_1-x_2)\,\mathrm dx_1\,\mathrm dx_2$ is the probability for
$X$ having a point in each of $B_1$ and $B_2$. For further details on
spatial point processes, see \citet{MW-04} and the references therein.

We view any vector $x=(x^{(1)},\ldots,x^{(d)})\in\mathbb R^d$ as a
column vector, and $\|x\|=\{(x^{(1)})^2+\ldots+(x^{(d)})^2\}^{1/2}$ as
its length. For ease of presentation, we assume that no pair of
distinct points $\{x_1,x_2\}\subset X$ is such that
$u=(x_1-x_2)/\|x_1-x_2\|$ is perpendicular to the $x^{(d)}$-axis. This
will happen with probability one for the models considered later in
this paper.

Let $\mathbb{S}^{d-1}=\{u=(u^{(1)},\ldots,u^{(d)})\in\mathbb R^d:
\,\|u\|=1\}$ be the unit sphere in $\mathbb R^d$ and
$e_d=(0,\ldots,0,1)$ its top point.  Denote $o=(0,\ldots,0)$ the
origin of $\mathbb R^d$.  Consider the $d$-dimensional cylinder with
midpoint $o$, radius $r>0$, height $2t>0$, and direction $e_d$:
\[
C(r,t)=\left\{x=(x^{(1)},\ldots,x^{(d)})\in
\mathbb R^d:\,(x^{(1)})^2+\cdots+(x^{(d-1)})^2 \le r^2,\,|x^{(d)}|\le t\right\}.
\] 
For $u\in\mathbb{S}^{d-1}$, denoting $\mathcal O_u$
an arbitrary $d\times d$ rotation matrix such that $u=\mathcal
O_{u}e_d$, then
\[
C_{u}(r,t)=\mathcal O_{u} C(r,t)
\]
is the $d$-dimensional cylinder with midpoint $o$, radius $r$,
height $2t$, and direction $u$.

\subsection{The cylindrical $K$-function}\label{s:cylK_def}

Recall that the second order reduced moment measure $\mathcal K$ with
structuring element $B\subset\mathbb R^d$, a bounded Borel set, is
given by
\[\mathcal K(B)=\int_Bg(x)\,\mathrm
dx\]
\citep[Section 4.1.2]{MW-04}. Ripley's $K$-function
\citep{ripley:76,ripley:77} is obtained when $B$ is a ball and it is not
informative about any kind of anisotropy in a spatial point pattern.

To detect preferred directions of linear structures in a spatial point
pattern, in particular a columnar structure, we propose a cylinder as
the structuring element and define the cylindrical $K$-function
  in the direction $u$ by
\begin{equation}\label{e:Kg}
K_{u}(r,t)=\int_{C_{u}(r,t)}g(x)\,\mathrm
dx,\quad u\in\mathbb{S}^{d-1},\quad r>0,\ t>0.
\end{equation}
Intuitively, $\rho K_{u}(r,t)$ is the mean number of further
points in $X$ within the cylinder with midpoint at the
typical point of $X$, radius $r$, and height $2t$ in the
direction $u$. For example, a stationary Poisson process is
isotropic, has $g=1$, and
\[
K_{u}(r,t)=2\omega_{d-1}r^{d-1}t,
\]
where $\omega_{d-1}=\pi^{(d-1)/2}/\Gamma\{(d+1)/2\}$ is the volume of
the $({d-1})$-dimensional unit ball. For $d=3$, $K_{(0,0,1)}$ is
similar to the space-time $K$-function in \citet{Diggleetal-95} and
\citet{GabrielDiggel-09}, when considering the $x^{(3)}$-axis as time and
the $(x^{(1)},x^{(2)})$-plane as space.

If $W\subset\mathbb R^d$ is an arbitrary Borel set with
$0<|W|<\infty$, then by standard methods \citep[Section 4.1.2]{MW-04}
\begin{equation}\label{e:ksum} 
  K_{u}(r,t)=\frac{1}{\rho^2|W|}\, E\left[\sum_{x_1,x_2\in
    X:\,x_1\not=x_2}\mathbb1\{x_1\in W,x_2- x_1\in C_{u}(r,t)\}\right],\quad r>0,\ t>0,
\end{equation}
where $\mathbb1$ denotes the indicator function, and by
stationarity the right-hand side does not depend on the choice of $W$.
This provides a more general definition of $K_{u}$, since
\eqref{e:ksum} does not require the existence of the pair correlation
function.  Equation \eqref{e:ksum} becomes useful when deriving
non-parametric estimates in Section~\ref{s:nonparest}.

\subsection{Non-parametric estimation}\label{s:nonparest}

Given a bounded observation window $W\subset\mathbb R^d$ and an observed point
pattern $\{x_1,\ldots,x_n\}\subset W$ with $n\ge2$ points,
 we consider non-parametric estimates of the form
\begin{equation}\label{e:Kest}
  \widehat K_{u}(r,t)=\frac{1}{\widehat{\rho^2}}\sum_{i\not= j} 
w_{u}(x_i,x_j)
  \mathbb1\{x_j-x_i\in C_{u}(r,t)\}.
\end{equation}
Here $\widehat{\rho^2}$ is a non-parametric estimate of $\rho^2$ and
$w_{u}$ is an edge correction factor. If $X$ is isotropic,
$K_{u}(r,t)$ does not depend on $u$ and this should affect the choice of
$\widehat K_{u}(r,t)$. On the other hand, as illustrated in the right
panel of Figure~\ref{fig:cylinder} and in Section~\ref{s:examplesK},
to detect a preferred direction of linearity in a spatial point
pattern, we suggest using an elongated cylinder with $t>r$ and
considering different directions $u$. Then we expect a largest value
of $\widehat K_{u}(r,t)$ to indicate the preferred direction, but a
careful choice of $r$ and $t$ may be crucial, cf.\
Section~\ref{s:extensions}. Furthermore, since $K_{u}=K_{-u}$, we
need only to consider the case where $u=(u^{(1)},\ldots,u^{(d)})$ is
on the upper unit-sphere, i.e., $u^{(d)}\ge0$.

Specifically, we use $\widehat{\rho^2}=n(n-1)/|W|^2$, see e.g.\
\citet{Illianetal-08}, and the translation correction factor
\citep{OhSt-81}
\begin{equation}\label{e:trans}
w_{u}(x_1,x_2) = 1/|W\cap W_{x_2-x_1}|
\end{equation}
where $W_{x}$ denotes translation of the set $W$ by a vector
$x\in\mathbb R^d$. Then, by Lemma~4.2 in \citet{MW-04}, if
$\widehat{\rho^2}$ is replaced by $\rho^2$ in \eqref{e:Kest}, we have
an unbiased estimate of $K_{u}$. As in
Figures~\ref{fig:cylinder}--\ref{fig:Data}, if $W$ is rectangular with
sides parallel to the axes and of lengths $a_1,\ldots,a_d>0$,
\[
|W\cap W_{x_2-x_1}| = \prod_{i=1}^d
\{a_i-|x_{2}^{(i)}-x_{1}^{(i)}|\},\quad x_1,x_2\in W,
\]
where $x_{j}^{(i)}$ denotes the $i$'th coordinate of $x_j$ $(j=1,2)$.

For $d=3$, $W=[0,a_1]\times[0,a_2]\times[0,a_3]$, and $u=(0,0,1)$,
another choice is a combined correction factor
\[w_{(0,0,1)}(x_1,x_2)=
\frac{1+\mathbb1\left(2x_{1}^{(3)}-x_{2}^{(3)}\not\in[0,a_3]\right)}
{a_3\left(a_1-\left|x_{2}^{(1)}-x_{1}^{(1)}\right|\right)
\left(a_2-\left|x_{2}^{(2)}-x_{1}^{(2)}\right|\right)}
\]
where the numerator is a temporal correction factor and the
denominator is the reciprocal of a spatial correction factor;
similarly we construct combined correction factors when $u=(1,0,0)$
or $u=(0,1,0)$. Instead of this spatial correction factor, which is
of a form similar to \eqref{e:trans}, \citet{Diggleetal-95} used an
isotropic correction factor, but this is only appropriate if $X$ is
isotropic in the $(x_1,x_2)$-plane. The temporal correction factor is
the same as that used in \citet{Diggleetal-95}.

We prefer the translation correction factor \eqref{e:trans}, since
this does not restrict the shape of $W$ and the choice of $u$.
In a simulation study with $d=3$, $W=[0,1]^3$, and $X$ a
Poisson line cluster point process as defined in
Section~\ref{s:PLCPP}, we obtained similar results when using the
translation and the combined correction factors.

\subsection{Examples}\label{s:examplesK}

\begin{figure}
\centering
\fbox{\includegraphics[height=4.5cm, width=4.6cm]{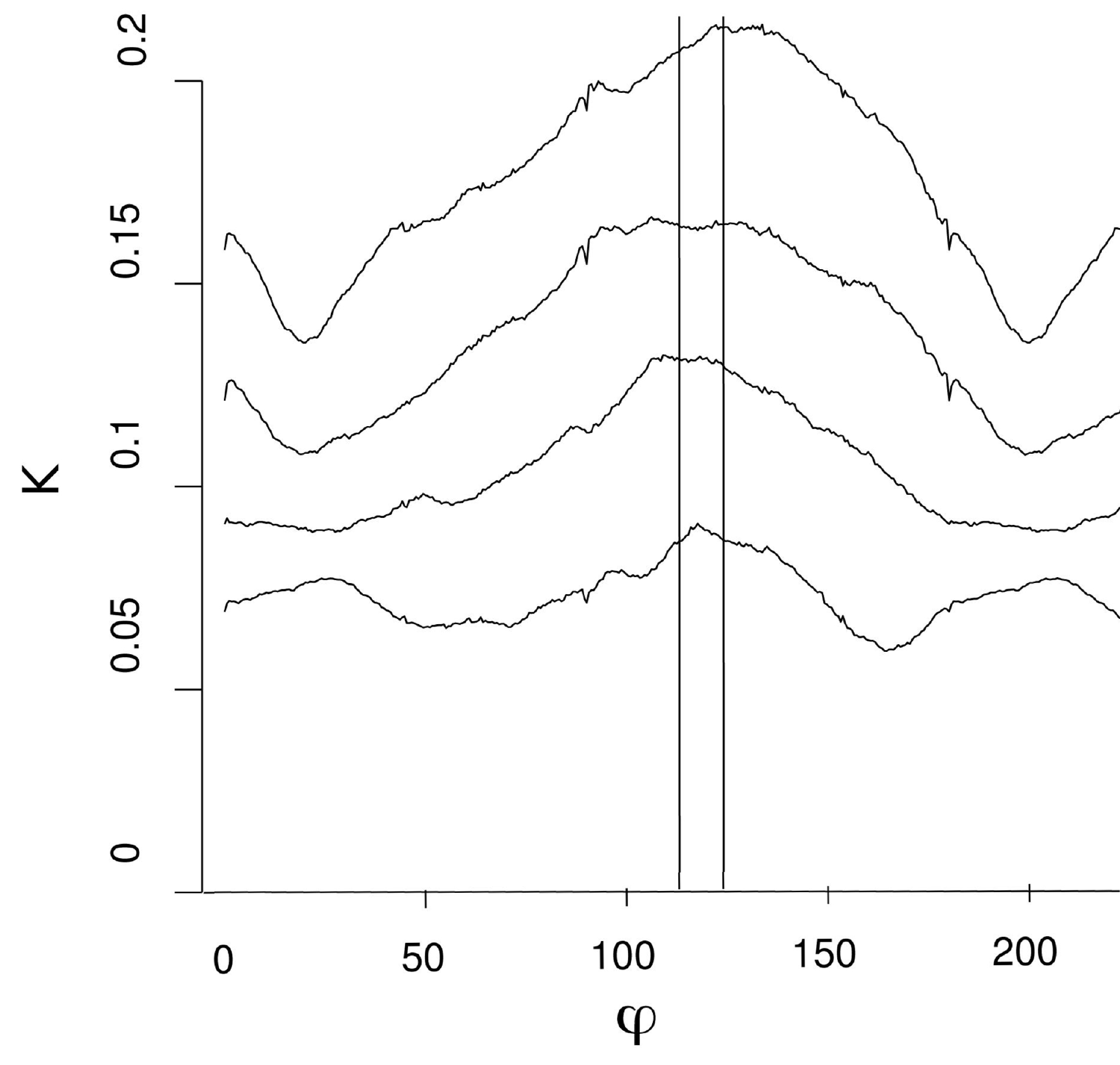}}
\fbox{\includegraphics[height=4.5cm, width=4.6cm]{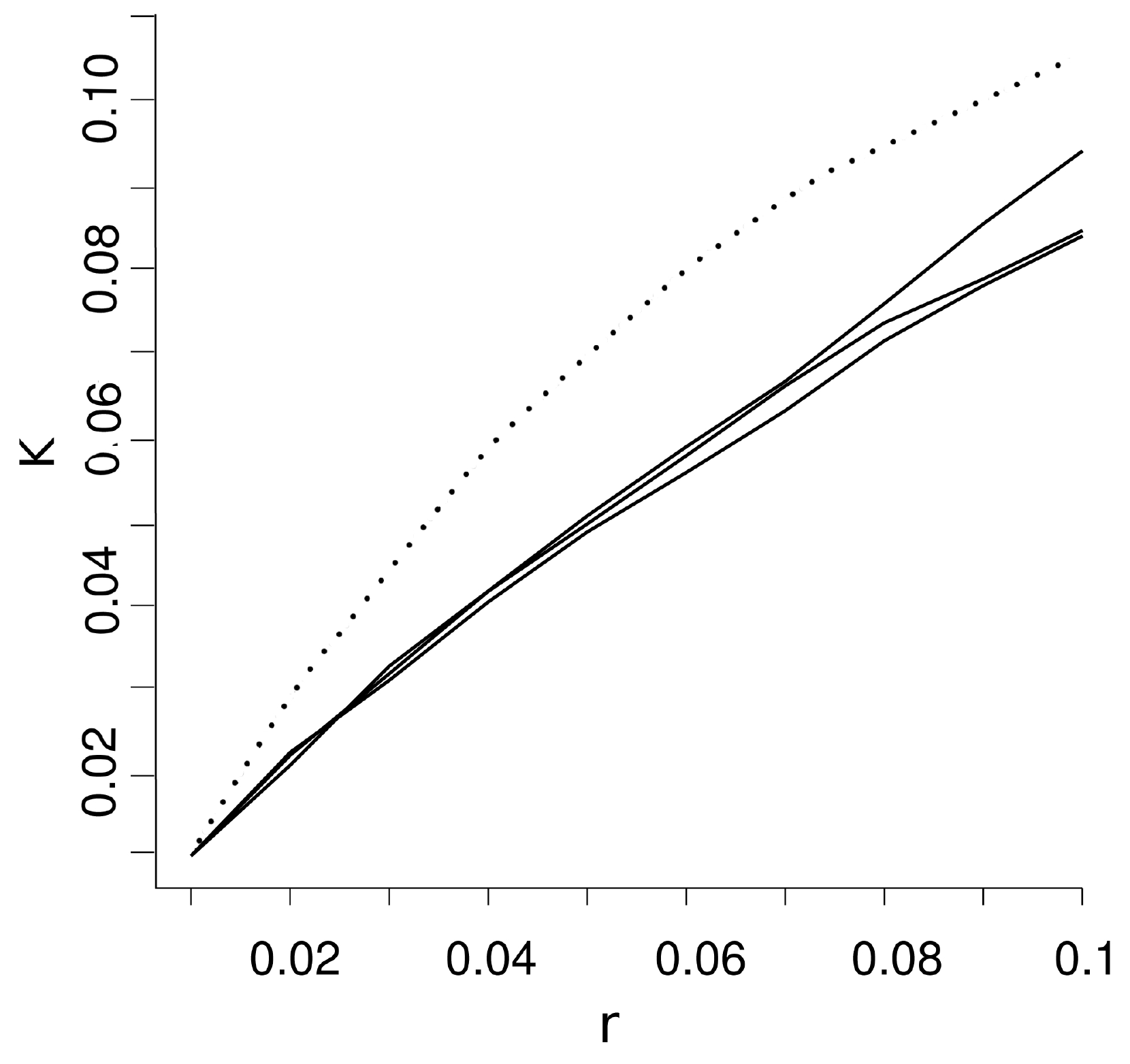}}
\fbox{\includegraphics[height=4.5cm, width=4.6cm]{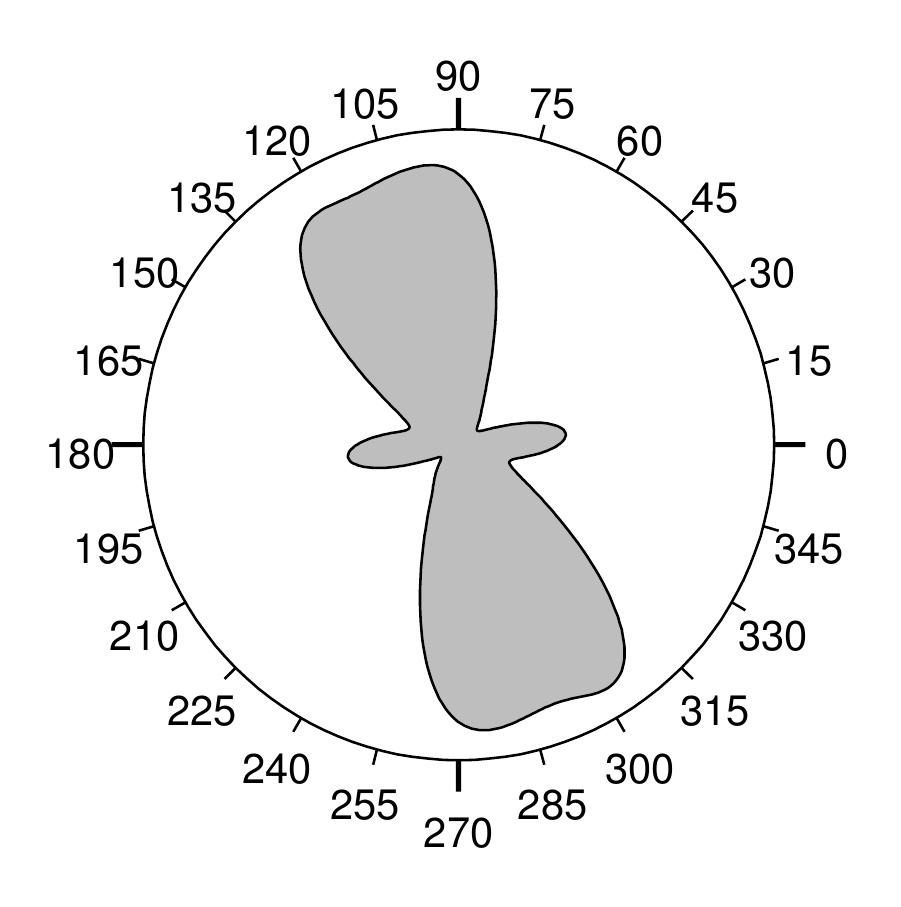}}
\caption{Chapel data set: The left panel shows the non-parametric
  estimate $\widehat K_{(\cos\varphi,\sin\varphi)}(r,t)$ versus
  $\varphi$ for four different combinations of $r$ and $t$ with the
  curves from the top to the bottom corresponding to $(r,t)=(0.2,0.4),
  (0.2,0.3), (0.1,0.3), (0.1,0.2)$ and the two vertical lines
  corresponding to $113^\circ$ and $124^\circ$. The middle panel shows
  $\widehat K_{(\cos\varphi,\sin\varphi)}(r,t)$ versus $r$ for
  different values of $\varphi$ and with $t=0.3$ with the solid curves
  from the top to the bottom corresponding to
  $20^\circ,45^\circ,170^\circ$, and the dotted curve corresponding to
  $\varphi=117^\circ$. The right panel shows a non-parametric estimate
  of the point pair orientation distribution function with $r_1=0.05$
  and $r_2=0.15$. For more details, see Section~\ref{s:examplesK}.}
\label{fig:aniso2D}
\end{figure}

\begin{figure}
\centering
\includegraphics[height=7cm, width=7cm]{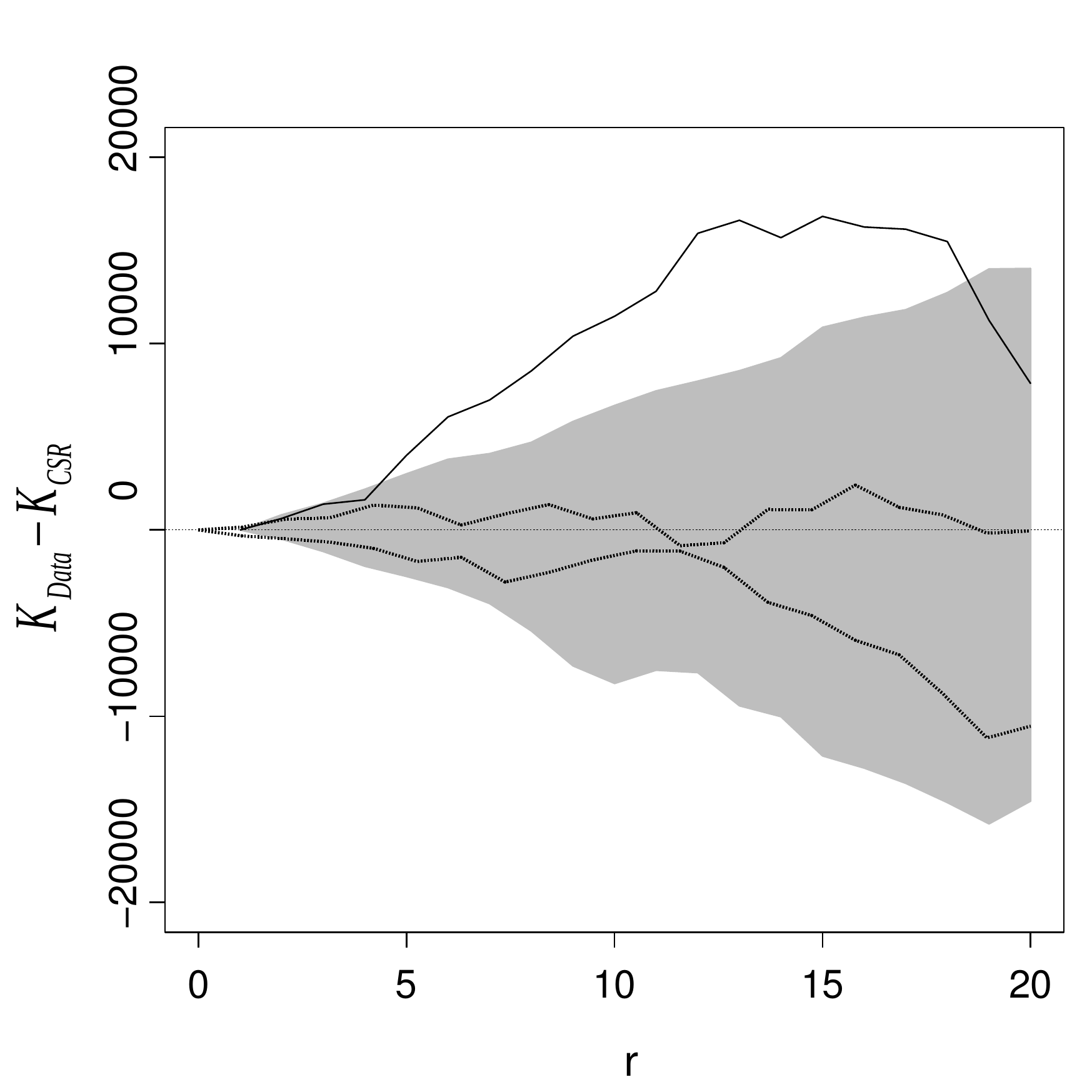}
\caption{Pyramidal cell data set: Non-parametric estimates $\widehat
  K_{u}(r,80)-160\pi r^2$ versus $r$ when $t=80$ and the cylinder is
  along the $x^{(1)}$-axis (dotted line), the $x^{(2)}$-axis (dotted
  line), or the $x^{(3)}$-axis (solid line). The grey region
  specifies a $95\%$ simultaneous rank envelope computed from 999
  simulations under complete spatial randomness. For more details, see
  Section~\ref{s:examplesK}.}
\label{fig:clyK_mini}
\end{figure}

Non-parametric estimates of the cylindrical $K$-function for the
two-dimensional chapel data set and the three-dimensional pyramidal
cell data set are shown in Figures~\ref{fig:aniso2D} and
\ref{fig:clyK_mini}, respectively. Below we comment on these
plots. Further examples are given in \citet{RSDRMN-15}.

To detect the main direction in the chapel point pattern, the left
panel in Figure~\ref{fig:aniso2D} shows, for four different
combinations of $r$ and $t$, i.e., $(0.1, 0.2)$, $(0.1, 0.3)$, $(0.2,
0.3)$ and $(0.2, 0.4)$, plots of $\widehat K_{u}(r,t)$ versus
$\varphi$, where $u = (\cos(\varphi), \sin(\varphi))$. These four
curves are approximately parallel, and a similar behaviour for other
choices of $r=0.05,0.1,0.15,0.2$ and $t=0.15,0.2,0.25,0.3,0.35,0.4$
with $t> 2r$ was observed.  In a previous analysis,
\citet{MoellerHokan-14} estimated the orientation of the chapel point
pattern to be between $113^{\circ}$ and $124^{\circ}$. This interval,
which is specified by the vertical lines in the left panel of
Figure~\ref{fig:aniso2D}, is in close agreement with the maximum of
the $\widehat K_{u}(r,t)$-curves. The middle panel in
Figure~\ref{fig:aniso2D} shows plots of $\widehat K_{u}(r,t)$ versus
$r$ when $t=0.3$. For the dotted curve in the middle panel,
$\varphi=117^{\circ}$ is the average of the four maximum points of
$\varphi$ corresponding to the four curves in the left panel, while
for the three other curves in the middle panel, values of $\varphi$
not included in the interval $[113^{\circ},124^{\circ}]$ have been
chosen. The clear difference between the dotted curve and the other
curves indicates a preferred direction in the point pattern which is
about $117^{\circ}$. This is also confirmed by the right panel in
Figure~\ref{fig:aniso2D}, which shows a non-parametric estimate of the
point pair orientation distribution function given by equation (14.53)
in \citet{StoyStoy-94} and implemented in {\texttt{spatstat}}
\citep{spatstat}. This is a kernel estimate which considers the
direction for each pair of observed points that lie more than
$r_1=0.05$ and less than $r_2=0.15$ units apart. For other values of
$r_2\le0.27$ we reached similar conclusions, but for higher values of
$r_2$ the pair orientation distribution function did not show a clear
preferred direction in the data.

By the minicolumn hypothesis, the pyramidal cell data set has a
columnar arrangement in the direction of the $x^{(3)}$-axis indicated
in Figure~\ref{fig:Data} \citep{RSDRMN-15}.
Figure~\ref{fig:clyK_mini} shows that the cylindrical $K$-function is
able to detect this kind of anisotropy: The three curves are $\widehat
K_{u}(r,80)-160\pi r^2$ for $0< r\le20$ and $u$ parallel to one of the
three main axes, where $160\pi r^2$ is the value of $K_{u}(r,80)$ under
complete spatial randomness, i.e.\ a stationary Poisson point
process model; we only make this comparison in order to see any
deviations from complete spatial randomness. The solid curve
corresponding to $u=(0,0,1)$, i.e., when the direction of the cylinder
is along the $x^{(3)}$-axis, is clearly different from the two other
cases where $u=(1,0,0)$ or $u=(0,1,0)$.  The grey region is a
so-called 95\% simultaneous rank envelope \citep{Myllymaki-16}
obtained from $999$ simulated realizations under complete spatial
randomness; \citet{Myllymaki-16} recommended 2499 simulations,
however, for the cylindrical $K$-function considered in
Figure~\ref{fig:clyK_mini}, 999 simulations seemed sufficient since
results were produced that were similar to those using 2499
simulations. Roughly speaking, under complete spatial randomness,
each of the estimated cylindrical $K$-functions is expected to be
within the grey region with estimated probability 95\%.  While the
curves for $u=(1,0,0)$ and $u=(0,1,0)$ are completely within the grey
region, the curve for $u=(0,0,1)$ is clearly outside for a large range
of $r$-values. In fact, for the null hypothesis of complete spatial
randomness, considering the rank envelope test \citep{Myllymaki-16}
based on $\widehat K_{u}(r,80)$ when $0<r\le20$ and $u=(0,0,1)$, the
$p$-value is estimated to be between $0.1\%$ and $0.18\%$, showing a
clear deviation from the null hypothesis of complete spatial
randomness.

These examples illustrate that the cylindrical $K$-function is a
useful functional summary statistic for detecting preferred directions
and columnar structures in a spatial point pattern. Particularly, for
the pyramidal cell data set and in accordance to the minicolumn
hypothesis, there is a pronounced columnar arrangement in the
direction of the $x^{(3)}$-axis.

\section{The Poisson line cluster point process}\label{s:PLCPP}

\subsection{Definition of Poisson line cluster point processes}\label{s:defprop}

Motivated by the analysis in Section~\ref{s:examplesK}, in particular
that of the pyramidal cell data set, we now introduce a model for a
point process $X$ with columnar structure. It is a
hierarchical construction, using various latent processes specified
briefly in (A1)--(C1) and later in more detail in (A2)--(C2), while
$X$ is given in (D1) and (D2) below; the left panel of
Figure~\ref{fig:cylinder} is helpful in this connection: It consists
of generating:
\begin{description}
\item[\textnormal{(A1)}] a Poisson line process $L=\{l_1,l_2,\ldots\}$
  of lines $l_i$, i.e., infinite, directed, and straight lines; only
  parts of these lines are shown in Figure~\ref{fig:cylinder};
\item[\textnormal{(B1)}] on each line $l_i$, a Poisson process $Q_i$ illustrated by unfilled points in Figure~\ref{fig:cylinder};
\item[\textnormal{(C1)}] a new point process $X_i$ obtained by random
  displacements in $\mathbb R^d$ of the points in $Q_i$ illustrated by
  filled points in Figure~\ref{fig:cylinder};
\item[\textnormal{(D1)}] finally, $X$ as the superposition of all the
  $X_i$.
\end{description}
Then we call $X$ a Poisson line cluster point process,
since its points cluster around the Poisson lines, and we call each
$X_i$ a cluster.




In connection to the more detailed conditions (A2)--(D2) we need the
following notation.
Let $\cdot$ denote the usual inner product on $\mathbb R^d$.  For
$u\in\mathbb{S}^{d-1}$, let $u^\bot=\{x\in\mathbb R^d:\, x\cdot u=0\}$
be the hyperplane perpendicular to $u$ and containing $o$,
$\lambda_{u^\bot}$ the $(d-1)$-dimensional Lebesgue measure on
${u^\bot}$, and $p_{u^\bot}(x)=x-(x\cdot u)u$ the orthogonal
projection of $x\in\mathbb R^d$ onto ${u^\bot}$.  Let $H={e_d^\bot}$
and $\lambda=\lambda_{e_d^\bot}$, i.e., $H$ is the hyperplane
perpendicular to the $x^{(d)}$-axis.  Further, let $k$ be a density
function with respect to Lebesgue measure on $\mathbb R^{d-1}$.  As in
Section~\ref{s:setting}, suppose we have specified for each
$u\in\mathbb{S}^{d-1}$ a $d\times d$ rotation matrix $\mathcal O_{u}$
such that $u=\mathcal O_{u}e_d$.  We then define a density with
respect to $\lambda_{u^\bot}$ by
\[
k_{u^\bot}\{\mathcal O_{u}(x^{(1)},\ldots,x^{(d-1)},0)\}
=k(x^{(1)},\ldots,x^{(d-1)}), \quad
(x^{(1)},\ldots,x^{(d-1)})\in\mathbb R^{d-1}.
\] 
In other words, when considering coordinates with respect to the $d-1$
first columns in $\mathcal O_{u}$, the distribution under $k_{u^\bot}$
is the same as under $k$.  Furthermore, for the line process $L$ we
use the so-called phase representation \citep{Stoyanetal-94} and
assume that with probability one, $L$ has no line contained in
$H$. Thereby a line $l=l(y,u)$ in $L$ corresponds to its direction
$u\in\mathbb{S}^{d-1}$ and its intersection point $y$ in $H$. Thus
$L=\{l_1,l_2,\ldots\}$ can be identified by a point process
$\Phi=\{(y_1,u_1),(y_2, u_2),\ldots\}\subset H\times\mathbb{S}^{d-1}$
such that $l_i=l(y_i,u_i)$ $(i=1,2,\ldots)$ and $\Phi\subset
H\times(\mathbb{S}^{d-1}\setminus H)$ almost surely.




In addition to (A1)--(D1) we assume that:
\begin{description}
\item[\textnormal{(A2)}] $\Phi$ is a Poisson process with intensity
  measure $\beta\lambda(\mathrm dy)M(\mathrm du)$, where $\beta>0$ is
  a parameter and $M$ is a probability measure on $\mathbb{S}^{d-1}$
  describing the direction of a typical line, where
  $M(\mathbb{S}^{d-1}\cap H)=0$; we assume that
  $\int1/|u^{(d)}|\,M(\mathrm du)<\infty$, where $u^{(d)}$ is the last
  coordinate of $u$; this assumption will be needed when we later in
  \eqref{e:L} specify the intensity and rose of directions of the line
  process;
\item[\textnormal{(B2)}] conditional on $\Phi$, we have that $
  Q_1,Q_2,\ldots$ are independent stationary Poisson processes
  on $l_1,l_2,\ldots$, respectively, with the same intensity
  $\alpha>0$;
\item[\textnormal{(C2)}] conditional on $\Phi$ and $Q_1,
  Q_2,\ldots$, we have that $X_1,X_2,\ldots$ are
  independent point processes and each $X_i$ is obtained by
  independent and identically distributed random displacements of the points in $Q_i$
  following the density $k_{u_i^\bot}$; thus $X_i$ is
  a Poisson process on $\mathbb R^d$ with intensity function
\begin{equation}\label{e:lambdai1}
  \Lambda_i(x)=
  \alpha k_{u_i^\bot}\{p_{u_i^\bot}(x-y_i)\},
  \quad x\in\mathbb R^d;
\end{equation}  
\item[\textnormal{(D2)}] hence the superposition $
  X=\cup_{i=1}^\infty X_i$ is a Cox process driven by
  $\Lambda=\sum_i\Lambda_i$, i.e.\ $X$ conditional on
  $\Phi$ is a Poisson process with intensity function
\begin{equation}\label{e:1}
  \Lambda(x)=\alpha
  \sum_{i=1}^\infty
  k_{u_i^\bot}\{p_{u_i^\bot}(x-y_i)\},
  \quad x\in\mathbb R^d.
\end{equation}
\end{description}
Some comments are in order.

The processes $L$, $\Lambda$, and $X$ are stationary, the
distribution of $L$ is given by $(\beta,M)$, and the
distribution of $X$ is determined by $(\beta,M,\alpha,k)$.

By (C2), conditional on $L$, for each line $l_i\in L$ and each point
$q_{ij}\in Q_i$, there is a corresponding point $x_{ij}\in X_i$
such that the random shift $z_{ij}=x_{ij}-q_{ij}$ follows the density
$k_{u_i}^\bot$.  We could have defined the Poisson line
cluster point process by letting the displacements follow a
distribution on $\mathbb{R}^d$ rather than a hyperplane, or more
precisely by letting $z_{ij}$ follow a density
\[
k_{u_i^\bot}\{p_{u_i^\bot}(z_{ij})\}
f_{u_i}\{z_{ij}-p_{u_i}(z_{ij})\}, \quad z_{ij}\in\mathbb R^d,
\] 
where $f_{u_i}$ is a density function with respect to Lebesgue measure
on the line $l_i-y_i=\{tu_i:t\in\mathbb R\}$.  However, since
the part of the displacements running along the line $l_i$ just
corresponds to independent displacements of a stationary Poisson
process, this will just result in a new stationary Poisson process
with the same intensity, see, e.g., Section~3.3.1 in \citet{MW-04},
and so there is essentially no difference.


The fact in (D2) that $X$ is a Cox process becomes important for the
calculations and the statistical methodology considered later in this
paper.


\subsection{Intensity and rose of directions for the Poisson line
  process}\label{s:int-rose}

We have specified the distribution of the Poisson line process $L$ by
$(\beta,M)$.  This is useful for computational reasons, but when
interpreting results it is usually more natural to consider the
intensity and the rose of directions of $L$, which we denote by
$\rho_L$ and $\mathcal R$, respectively. Formal definitions of these
concepts are given in Appendix~A, where it is shown
that for any Borel set $B\subseteq \mathbb S^{d-1}$,
\begin{equation}\label{e:L}
\rho_L=\beta\int1/|u^{(d)}|\,M(\mathrm du),\quad 
\mathcal R(B)=\int_B1/|u^{(d)}|\,M(\mathrm du)\bigg/
\int1/|u^{(d)}|\,M(\mathrm du).
\end{equation} 
In words, $\rho_L$ is the mean length of lines in $L$ within
any region of unit volume in $\mathbb R^d$, and $\mathcal R$ is the
distribution of the direction of a typical line in $L$, see,
e.g., \citet{Stoyanetal-94}.

Equation \eqref{e:L} establishes a one-to-one correspondence between
$(\rho_L,\mathcal R)$ and $(\beta,M)$, where
\begin{equation}\label{e:Lstar}
\beta=\rho_L\int|u^{(d)}|\,\mathcal R(\mathrm du),\quad
M(B)=\int_B|u^{(d)}|\,\mathcal R(\mathrm du)\bigg/
\int|u^{(d)}|\,\mathcal R(\mathrm du).
\end{equation}
Consequently, we can choose $\rho_L$ as any positive and finite
parameter, and $\mathcal R$ as any probability measure on $\mathbb
S^{d-1}$.  Moreover, $\beta\le\rho_L$ where the equality only holds
when $\mathcal R$ is concentrated with probability one at
$\pm{e_d}$. We call this special case the degenerate Poisson line
cluster point process.

For the rose of directions, we later use a von Mises--Fisher
distribution with concentration parameter $\kappa\ge0$ and mean
direction $\mu\in\mathbb{S}^{d-1}$.  This has a density
$f(\cdot\mid\mu,\kappa)$ with respect to the surface measure
on $\mathbb{S}^{d-1}$:
\begin{equation}\label{e:hdensity}
  f(u\mid\mu,\kappa)=c_{d}(\kappa)\exp(\kappa\mu\cdot
  u),\quad c_{d}(\kappa)=\frac {\kappa^{d/2-1}} {(2\pi)^{d/2}I_{d/2-1}(\kappa)},
  \quad u\in\mathbb{S}^{d-1},
\end{equation}
where $I_{d}$ denotes the modified Bessel function of the first kind
and order $d$.  Note that $L$ and $X$ are then isotropic if and only
if $\kappa=0$, in which case the choice of $\mu$ plays no
role.  For $\kappa>0$, the directions of the lines in $L$ are
concentrated around $\mu$, and so the clusters in $X$ have
preferred direction $\mu$. When $\mu=\pm{e_d}$,
in the limit as $\kappa\rightarrow\infty$, we obtain the degenerate
Poisson line cluster point process.


\subsection{Finite versions of the Poisson line cluster point process and
  simulation}\label{s:simulation}

Suppose we want to simulate the Poisson line cluster point process
within a bounded region $W\subset\mathbb R^d$, i.e.\ the restriction
$X_W=X\cap W$. Then we need a finite approximation of
$\Phi$, which will also be used when we later discuss
Bayesian inference, as follows. Consider a bounded region
$W_{\text{ext}}\supseteq W$ and let
\[
S=\{(y,u)\in H\times\mathbb S^{d-1}: l(y,u)\cap
W_{\text{ext}}\not=\emptyset\}
\] 
be the set of all lines hitting $W_{\text{ext}}$. We want to choose
$W_{\text{ext}}$ as small as possible but so that it is very unlikely
that for some line $l_i\in L$ with $(y_i,u_i)\not\in S$, $X_i$ has a
point in $W$. Then our finite approximation is $\Phi_S=\Phi\cap S$,
and (i) we simulate $\Phi_S$, and (ii) conditional on $\Phi_S$ we make
an approximate simulation of $X_W$ as a Poisson process with intensity
function
\begin{equation}\label{e:vedikke}
\Lambda_W(x)=\alpha
\sum_{(y,u)\in\Phi_S}
k_{u^\bot}\{p_{u^\bot}(x-y)\},
\quad x\in W,
\end{equation}
cf.\ \eqref{e:1}. We detail (i)--(ii) below.

Here (ii) is rather straightforward: Suppose we have simulated
$\Phi_S$ and consider any $(y_i,u_i)\in\Phi_S$.  The projection of $W$
onto $l_i$ is the bounded set $l_{W,i}=\{x\in l_i: (x+u_i^\bot)\cap
W\not=\emptyset\}$.  In accordance to (B2), we simulate a Poisson
process $Y_{W,i}$ with intensity $\alpha$ on $l_{W,i}$.  Displacing
the points in $Y_{W,i}$ as described in (C2) we obtain a Poisson
process $X_{W,i}$ with intensity function \eqref{e:lambdai1} but
restricted to $\cup_{x\in l_{W,i}}(x+u^\bot)$.  The approximate
simulation of $X_W$ is then given by
$\cup_{(y_i,u_i)\in\Phi_S}X_{W,i}\cap W$.
 
In (i) we assume for simplicity and specificity that $\mathcal R$
follows the von Mises--Fisher density \eqref{e:hdensity}.  Denote
$\nu_{d-1}$ the surface measure on $\mathbb S^{d-1}$. Then \eqref{e:Lstar} implies that 
\[
\beta\lambda(\mathrm dy)\mu(\mathrm d\mu)=\rho_Lf(y\mid\mu,\kappa)|u^{(d)}|\lambda(\mathrm dy)\nu_{d-1}(\mathrm du),
\]
i.e., $\Phi_S$ is a Poisson process on $S$ with intensity function
\[\chi(y,u\mid\rho_L,\mu,\kappa)=
\rho_L|u^{(d)}|f(u\mid\mu,\kappa)\] with respect to the
measure $\lambda(\mathrm dy)\nu_{d-1}(\mathrm du)$.
First, we therefore simulate the Poisson distributed counts
$\#\Phi_S$ with mean $\rho_LI(\mu,\kappa)$ where
\begin{equation*}
  I(\mu,\kappa)=\int|u^{(d)}|
  f(u\mid\mu,\kappa)\mathrm d
  \lambda(\mathrm dy)\nu_{d-1}(\mathrm du)=
  \int\lambda(J_{u})f(u\mid\mu,\kappa)\nu_{d-1}(\mathrm du)
\end{equation*}
where $J_{u}=\{y\in H:l(y,u) \cap W_{\text{ext}}\not=\emptyset\}$.
Second, we simulate each $(y,u)\in\Phi_S$ with density proportional to
$|u^{(d)}|f(u\mid\mu,\kappa)$ for $y\in J_{u}$ and
zero otherwise. Here we use rejection sampling.

For example, if $d=2$ and $W_{\text{ext}}=[-a,a]^2$ is a square
centered at the origin, then for $u=(\cos\varphi,\sin\varphi)$, we
have $J_{u}=J_{\varphi}\times\{0\}$ with
\begin{equation}\label{e:J}
J_{\varphi}=\left\{\begin{array}{ll}
\left[-a\cot\varphi-a,a\cot\varphi+a\right], & \mbox{$0<\varphi\le\pi/2$ or $\pi<\varphi\le 3\pi/2$},\\
\left[a\cot\varphi-a,a-a\cot\varphi\right],  & \mbox{$\pi/2\le\varphi<\pi$ or $3\pi/2\le\varphi<2\pi$}.
\end{array}
\right. 
\end{equation}
Further,
\begin{equation}\label{e:lambda_J}
\lambda(J_{u})=\left\{\begin{array}{ll}
2a+2a\cot\varphi, & \mbox{$0<\varphi\le\pi/2$ or $\pi<\varphi\le 3\pi/2$},\\
2a-2a\cot\varphi, & \mbox{$\pi/2\le\varphi<\pi$ or $3\pi/2\le\varphi<2\pi$},
\end{array}
\right. 
\end{equation}
and 
\begin{align}\label{e:I}
I(\mu,\kappa)=&
2a\int_0^{\pi/2}(\sin\varphi+\cos\varphi)
f(u\mid\mu,\kappa)\,\mathrm d\varphi+
2a\int_{\pi/2}^\pi(\sin\varphi-\cos\varphi)
f(u\mid\mu,\kappa)\,\mathrm d\varphi\notag\\
&-2a\int_\pi^{3\pi/2}(\sin\varphi+\cos\varphi)
f(u\mid\mu,\kappa)\,\mathrm d\varphi-
2a\int_{3\pi/2}^{2\pi}(\sin\varphi-\cos\varphi)
f(u\mid\mu,\kappa)\,\mathrm d\varphi,
\end{align}
which can be evaluated by numerical methods. Furthermore, for
$\mu=(\cos\theta,\sin\theta)$ and $y=(y^{(1)},y^{(2)})$,
the unnormalized density $|u^{(d)}|f(u\mid\mu,\kappa)=
\mathbb1(y^{(1)}\in
J_{\varphi})|\sin\varphi|\exp\{\kappa\cos(\varphi-\theta)\}$ is just
with respect to Lebesgue measure $\mathrm dy^{(1)}\,\mathrm d\varphi$ on
$\mathbb R\times[0,2\pi)$.  Finally, when doing rejection sampling, we
propose $\varphi$ from $f(\cdot\mid\mu,\kappa)$ and $y^{(1)}$ from
the uniform distribution on $J_{\varphi}$, and accept $(\varphi,y^{(1)})$
with probability $|\sin\varphi|$.

\subsection{Moments of the Poisson line cluster point process}\label{s:mom}

Since $X$ is a Cox process with driving random intensity $\Lambda$,
moment properties of $X$ are determined by moment properties of
$\Lambda$. This section focuses on first and second order moments.

The Poisson line cluster point process $X$ has intensity $\rho$ and
pair correlation function $g$
\begin{equation}\label{e:1a}
\rho = E\{\Lambda(o)\},\quad
\rho^2g(x) = E\{\Lambda(o)\Lambda(x)\},\quad
x\in\mathbb R^d. 
\end{equation}
Appendix~B verifies that 
\begin{equation}\label{e:2}
\rho=\alpha\rho_L 
\end{equation}
and
\begin{equation}\label{e:4}
g(x)=1+\frac{1}{\rho_L}\int 
k_{u^\bot}*\tilde k_{u^\bot}\{p_{u^\bot}(x)\}
\,\mathcal R(\mathrm du),
\quad x\in\mathbb R^d,
\end{equation}
where $\tilde k_{u^\bot}\{p_{u^\bot}(x)\}
=k_{u^\bot}\{-p_{u^\bot}(x)\}$ 
and $*$ denotes convolution, i.e.,
\[k_{u^\bot}*\tilde k_{u^\bot}\{p_{u^\bot}(x)\}= \int k_{u^\bot}\{p_{
  u^\bot}(x)-y\} \tilde k_{u^\bot}(y) \,\lambda_{u^\bot}(\mathrm
dy).\] Thus $g>1$, reflecting the clustering of the Poisson line
cluster point process. Evaluation of the integral in \eqref{e:4} may
require numerical methods.  For example, if
$k(\cdot)=f(\cdot\mid\sigma^2)$ is the density of the
$(d-1)$-dimensional zero-mean isotropic normal distribution with
variance $\sigma^2>0$, then
\begin{equation}\label{e:conv}
k_{u^\bot}*\tilde k_{u^\bot}\{p_{u^\bot}(x)\}=
\exp\left\{-\|p_{u^\bot}(x)\|^2/
\left(4\sigma^2\right)\right\}/\left(4\pi\sigma^2\right)^{(d-1)/2}.
\end{equation} 

\subsection{Moment based inference}\label{s:infmoments}

The likelihood for a parametric Poisson line cluster point process
model is complicated because of the hidden line process and the hidden
point processes on the lines, though it can be approximated using a
missing data Markov chain Monte Carlo approach, see, e.g.,
\citet{MW-04}. A Bayesian Markov chain Monte Carlo approach is used
in Section~\ref{s:Bayesinf} where the missing data is included into
the posterior.  Simpler procedures for parameter estimation are
composite likelihood \citep{Guan-06,MW-07} and minimum contrast
methods \citep{mincontrast-84} based on \eqref{e:2}-\eqref{e:4}. Since
$g$ is hard to compute in general, this section focuses on such simple
procedures in the special case of a degenerate Poisson line cluster
point process which e.g. is a relevant model for the pyramidal cell
data set shown in Figure~\ref{fig:Data}. For specificity we assume as in
\eqref{e:conv} that $k(\cdot)=f(\cdot\mid\sigma^2)$.  Then the unknown
parameters are $\beta=\rho_L>0$, $\alpha>0$, and $\sigma^2>0$.


Suppose a realization of $X_W$ is observed within a region of the
product form $W=D\times I$, where $D\subset\mathbb R^{d-1}$ and
$I\subset\mathbb R$ are bounded sets.  To estimate the unknown
parameters we notice the following. Let $X_I$ denote the projection of
$X\cap(\mathbb R^{d-1}\times I)$ onto $H$. Since we consider a
degenerate Poisson line cluster point process, the $x^{(d)}$-coordinates
of the points in $X_W$ are independent and identically distributed
uniform points on $I$ which are independent of $X_I$. Thus $X_I$ is a
sufficient statistic for $(\rho_L,\alpha,\sigma^2)$.  Note that $X_I$
is a Cox process driven by the random intensity function
\[
\Gamma(x)=\alpha|I|\sum_{i=1}^\infty f(x-y_i\mid\sigma^2)
,\quad x\in H,
\] 
where $\Phi=\{(y_1,e_d),(y_2,e_d),\ldots\}$ can be identified by the
stationary Poisson process $\{y_1,y_2,\ldots\}$ on $H$
with intensity $\rho_L$. Therefore $X_I$ is a modified Thomas
process \citep{MW-04} with intensity $\rho_I=\rho|I|$ and by
\eqref{e:4}-\eqref{e:conv} pair correlation function
\[
g_I(x)=1+\frac{1}{\left(4\pi\sigma^2\right)^{(d-1)/2}\rho_L}
\exp\left(-\frac{\|x\|^2}{4\sigma^2}\right),\quad
x\in H.
\]

Parameter estimation based on $(\rho_I,g_I)$ and using a composite
likelihood or a minimum contrast method is straightforward
\citep{MW-07}. Then, when checking a fitted Thomas process, we should
not reuse the intensity and the pair correlation function. Below we
use instead the functional summary statistics the empty space
function $F$, the nearest-neighbour function $G$, and the $J$-function
\citep{MW-04}.



\begin{figure}
\centering 
\includegraphics[height=5.8cm, width=5.7cm]{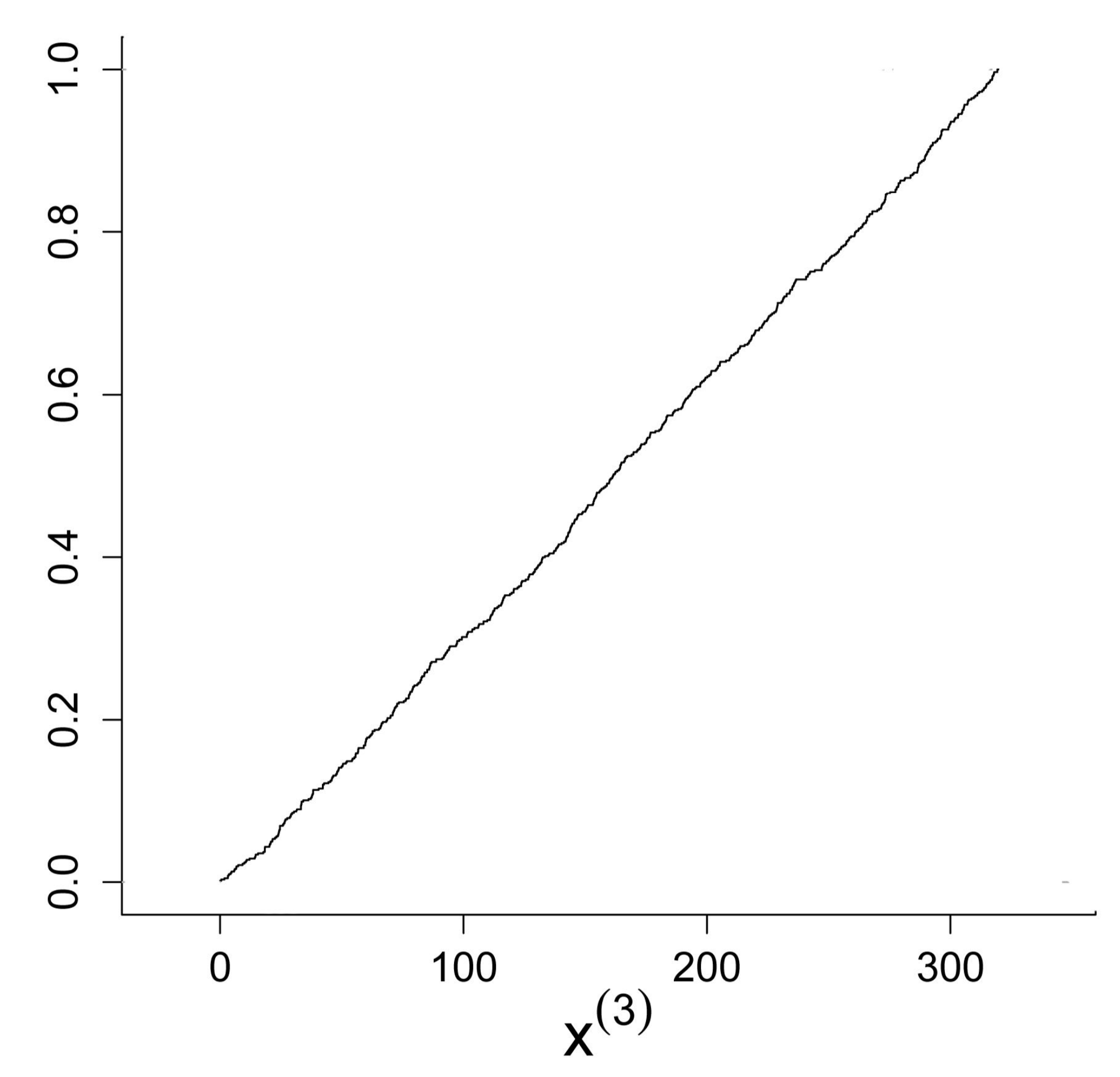}
\includegraphics[height=6cm, width=6cm]{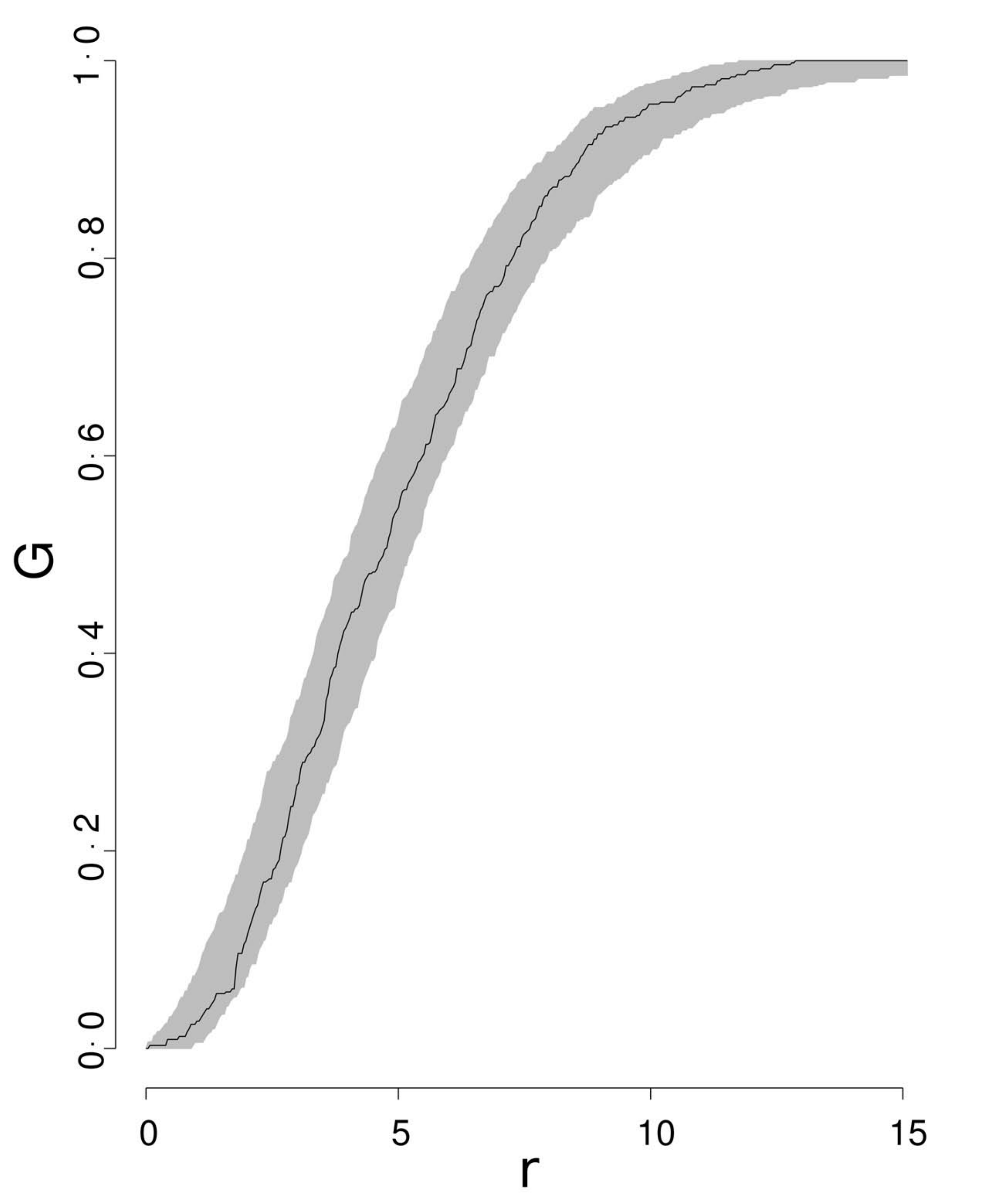}\\
\includegraphics[height=6cm, width=6cm]{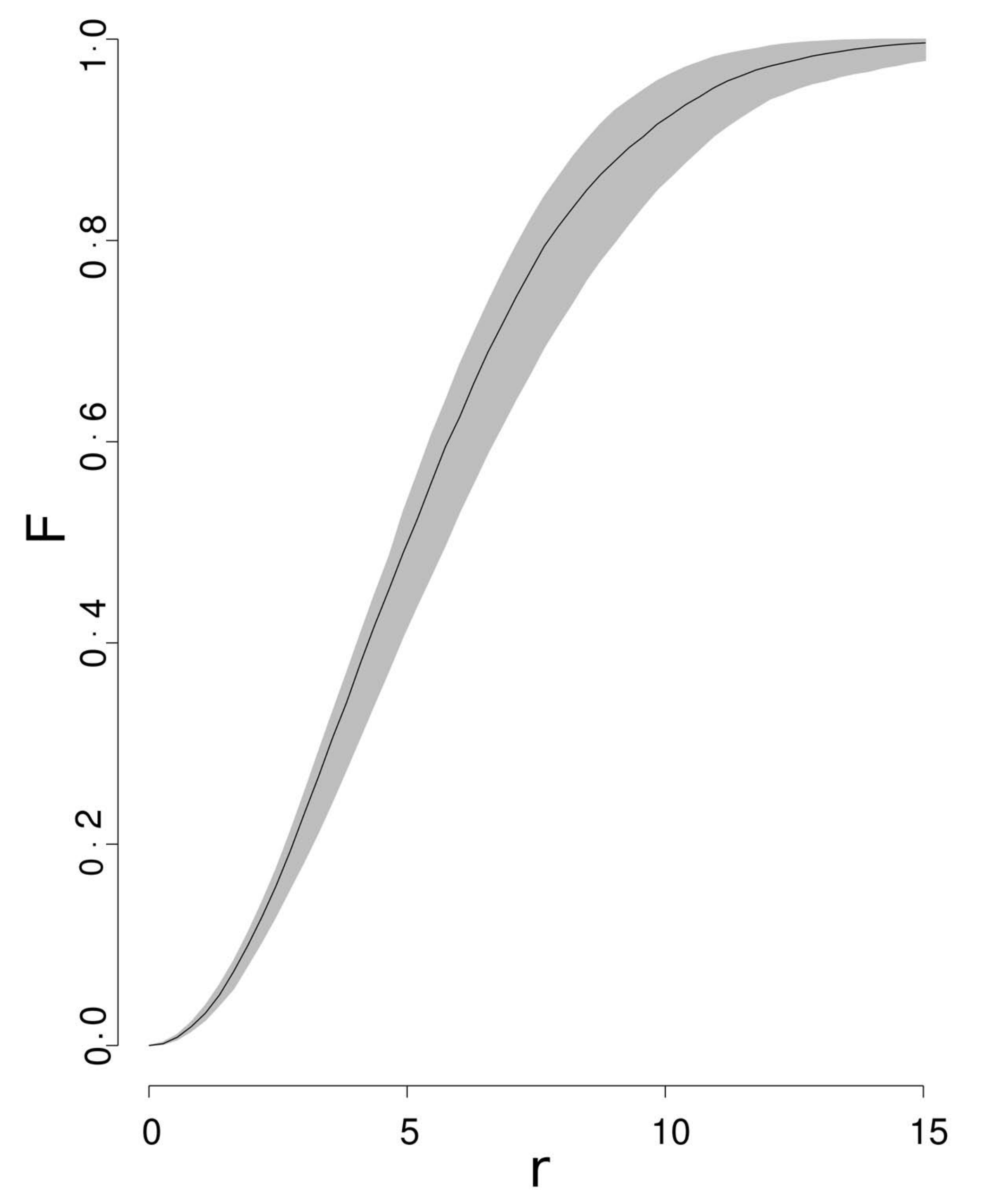}
\includegraphics[height=6cm, width=6cm]{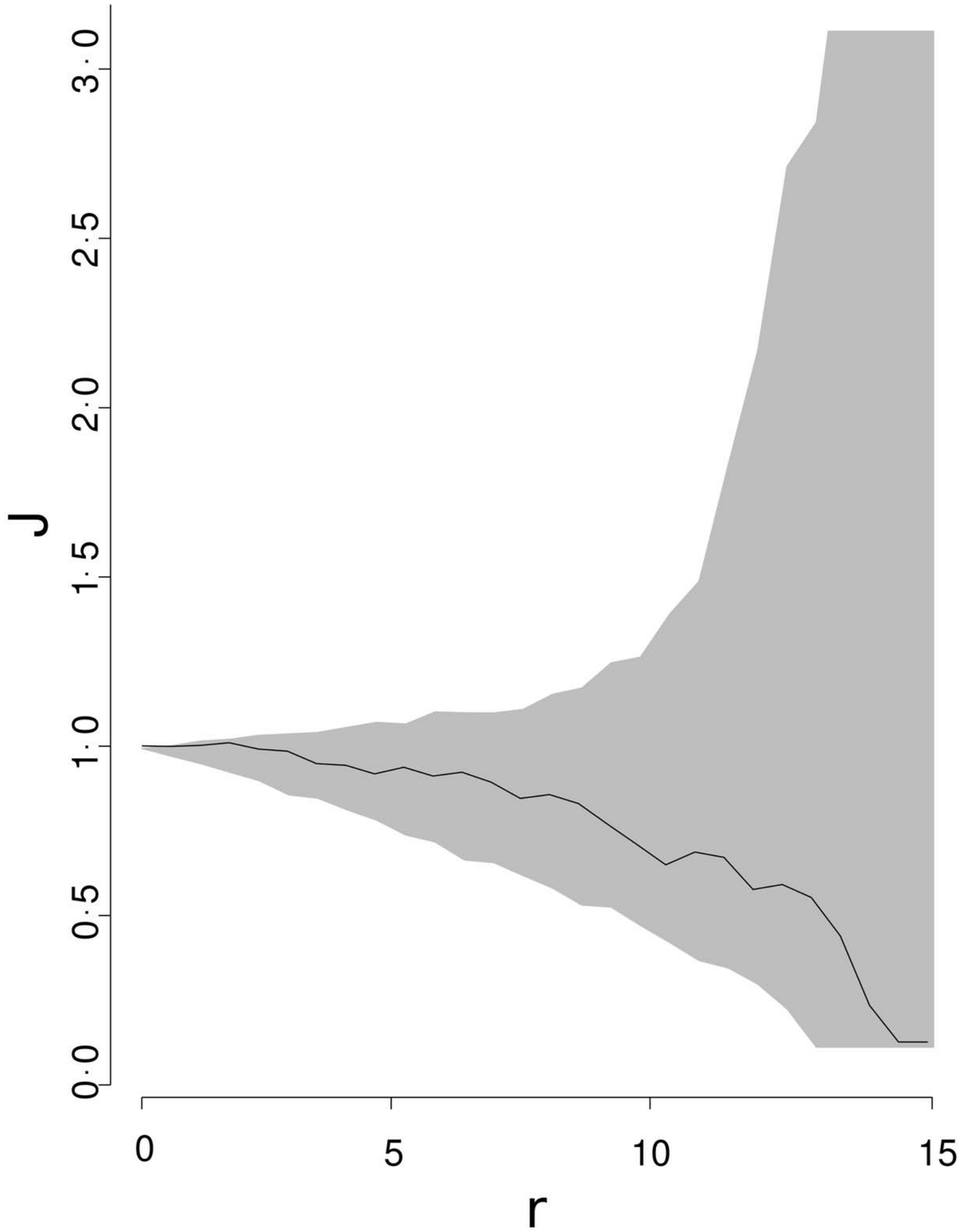}
\caption{Summary statistics for the pyramidal cell point pattern data
  set: Empirical cumulative distribution function of the
  $x^{(3)}$-coordinates of the pyramidal cell point pattern data set
  (top left), and non-parametric estimates of $G$ (top right), $F$
  (bottom left), and $J$ (bottom right) for the projected pyramidal
  cell point pattern onto $D$ (solid lines), together with $95\%$
  simultaneous rank envelopes (gray regions) computed from 4999
  simulated point patterns under the fitted Thomas process.}
\label{fig:fit}
\end{figure}

As an example, for the three-dimensional pyramidal cell data set in
Figure~\ref{fig:Data} in accordance with the minicolumn hypothesis, we
consider a degenerate Poisson line cluster point process. This has a
columnar arrangement in the direction of the $x^{(3)}$-axis and the
observation window is of the same form as described above with
$D=[0,508]\times[0,138]$ and $I=[0,320]$.  The first panel in
Figure~\ref{fig:fit} shows the empirical cumulative distribution
function of the $x^{(3)}$-coordinates of the pyramidal cell point
pattern data set; there is no clear indication of a deviation from a
uniform distribution, in agreement with our stationarity assumption.

When fitting the modified Thomas process for the projected point
pattern onto $D$, for both composite likelihood and minimum contrast
estimation, we used the {\texttt{spatstat}} \citep{spatstat} function
{\texttt{kppm}}.  We obtained the minimum contrast estimates
$\hat\rho_L=0.024$, $\hat\alpha=0.37/320=0.0012$, and
$\hat\sigma^2=15.04$, and similar estimates were obtained by a
composite likelihood method. The three last panels in
Figure~\ref{fig:fit} show the non-parametricly estimated $F$, $G$, and
$J$-functions for the projected pyramidal cell point pattern onto $D$,
together with $95\%$ simultaneous rank envelopes computed from 4999
simulated point patterns under the fitted Thomas process. The
$p$-values for the rank envelope test \citep{Myllymaki-16} for $G$,
$F$, and $J$-functions are within the intervals $[0.851, 0.852]$,
$[0.732, 0.733]$, and $[0.623, 0.625]$, respectively, providing no
evidence against the fitted model.

\subsection{Bayesian inference}\label{s:Bayesinf}

Suppose a non-empty realization $X_W=\{x_1,\ldots,x_n\}$ is our data,
where $W\subset\mathbb R^d$ is a bounded observation window, and we
model $X$ as a Poisson line cluster point process with
$k_{u^\bot}(y)=f(y\mid\sigma^2)$ and $\mathcal R$ following the von
Mises--Fisher density $f(u\mid\mu,\kappa)$ given by
\eqref{e:hdensity}.  As in Section~\ref{s:simulation} we need a finite
representation $\Phi_S$ of $\Phi$ which we treat as a latent process.
This section considers a Bayesian Markov chain Monte Carlo missing
data approach for the missing data $\Phi_S$ and the unknown parameters
$\rho_L>0$, $\mu\in\mathbb S^{d-1}$, $\kappa>0$, $\alpha>0$, and
$\sigma^2>0$ .


Imagining that also a realization $\Phi_S=
\{(y_1,u_1),\ldots,(y_k,u_k)\}$ had been observed, we detail below the
calculation of the likelihood $l[\rho_L,\mu,\kappa,\alpha,\sigma^2\mid
\{x_1,\ldots,x_n\},\{(y_1,u_1),\ldots, (y_k,u_k)\}]$. For the
parameters, we assume independent prior densities
$p(\rho_L),p(\mu),p(\kappa),p(\alpha),p(\sigma^2)$; further prior
specifications are given below. Hence the posterior density is
\begin{align}
&p[\rho_L,\mu,\kappa,\alpha,\sigma^2,
\{(y_1,u_1),\ldots,
(y_k,u_k)\}\mid
\{x_1,\ldots,x_n\}]\nonumber\\
&\propto
l[\rho_L,\mu,\kappa,\alpha,\sigma^2\mid
\{x_1,\ldots,x_n\},\{(y_1,u_1),\ldots,
(y_k,u_k)\}]
p(\rho_L)p(\mu)p(\kappa)p(\alpha)p(\sigma^2).
\label{e:post}
\end{align}

As a first ingredient of the likelihood, using the approximation
$\Phi_S$ of $\Phi$, we also approximate $X_{W}$
by a finite Cox process $X_{W,S}$ with driving random
intensity function $\Lambda_W$ given by \eqref{e:vedikke} with
$k_{u^\bot}(\cdot)=f(\cdot\mid\sigma^2)$. Conditional on
$\Phi_S$, $X_{W,S}$ is absolutely continuous with
respect to the unit rate Poisson process on $W$, with density
\begin{equation}\label{e:den1}
f[\{x_1,\ldots,x_n\}\mid\Phi_S,\alpha,\sigma^2]=
\exp\left\{|W|-
\int_W\Lambda_W(x\mid\Phi_S,\alpha,\sigma^2)
\,\mathrm dx\right\}
\prod_{i=1}^n\Lambda_W(x_i\mid\Phi_S,\alpha,\sigma^2) 
\end{equation}
for finite point configurations $\{x_1,\ldots, x_n\}\subset W$.

For the second ingredient of the likelihood, notice that the
distribution of $\Phi_S$ is absolutely continuous with respect
to the distribution of a natural reference process $\Phi_{0,S}$
defined as the Poisson process on $S$ with intensity function
$\chi_0(y,u)=|u^{(d)}|\Gamma(d/2)/(2\pi^{d/2})$ with
respect to the measure $\lambda(\mathrm dy)\nu_{d-1}(\mathrm
du)$, cf.\ Section~\ref{s:simulation}. This reference process
corresponds to the case of an isotropic Poisson line process with unit
intensity.  The density of $\Phi_S$ with respect to the
distribution of $\Phi_{0,S}$ is
\begin{align*}
&f[\{(y_1,u_1),\ldots,
(y_k,u_k)\}\mid\rho_L,\mu,\kappa]\\
&=\exp\left[\int_S\left\{\chi_0(y,u)-\chi(y,
u\mid\rho_L,\mu,\kappa)\right\}\lambda(\mathrm dy)\nu_{d-1}(\mathrm du)\right]
\prod_{j=1}^k \frac{\chi(y_j,u_j\mid\rho_L,\mu,\kappa)}
{\chi_0(y_j,u_j)}
\end{align*}
for finite point configurations $\{(y_1,u_1),\ldots,(y_k,u_k)\}\subset
S$.  That is, using the notation in Section~\ref{s:simulation},
\begin{align}
&f[\{(y_1,u_1),\ldots,(y_k,u_k)\}
\mid\rho_L,\mu,\kappa]\nonumber\\
&\propto\exp\left\{-\rho_LI(\mu,\kappa)\right\}
\prod_{j=1}^k \left\{\frac{2\pi^{d/2}}{\Gamma(d/2)}
\rho_L f(u_j\mid\mu,\kappa)
\mathbb1(y_j\in J_{u_j})\right\},
\label{e:den2}
\end{align}
where we have omitted a constant not depending on the parameters. 

Combining \eqref{e:den1}--\eqref{e:den2} we obtain the approximate
likelihood
\begin{align}
l[\rho_L,&\mu,\kappa,\alpha,\sigma^2\mid
\{x_1,\ldots,x_n\},\{(y_1,u_1),\ldots,
(y_k,u_k)\}]\nonumber\\
=\,&
\exp\left\{|W|-
\int_W\Lambda_W(x\mid\Phi_S,\alpha,\sigma^2)
\,\mathrm dx\right\}
\prod_{i=1}^n\Lambda_W(x_i\mid\Phi_S,\alpha,\sigma^2)
\nonumber\\
&
\times\exp\left\{-\rho_LI(\mu,\kappa)\right\}
\prod_{j=1}^k \left\{\frac{2\pi^{d/2}}{\Gamma(d/2)}
\rho_L f(u_j\mid\mu,\kappa)
\mathbb1(y_j\in J_{u_j})\right\}.
\label{(20)}
\end{align} 
Inserting this into \eqref{e:post}, we notice that the posterior
density is analytically intractable.  A hybrid Markov chain Monte
Carlo algorithm or Metropolis within Gibbs algorithm, see e.g.\
\citet{Gilksetal-96}, for posterior simulations is proposed in
Appendix~C. Briefly, the algorithm alternates
between updating each of the parameters and the line process, using a
birth-death-move Metropolis Hastings algorithm for the line process.

\begin{figure}
\centering
\includegraphics[width=\textwidth]{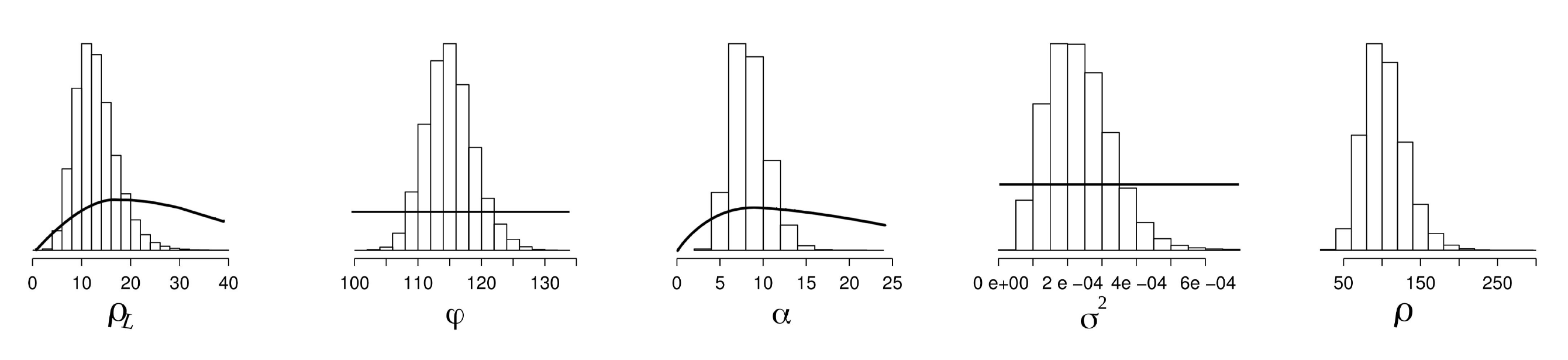}
\caption{Prior and posterior distributions: The first four panels show
  the unnormalized prior density (solid line) and histogram for the
  posterior distribution of $\rho_L$, $\varphi$, $\alpha$, and
  $\sigma^2$, respectively. The final panel shows the histogram for
  the posterior distribution of $\rho$.}
\label{fig:marg_param_hist}
\end{figure}

To illustrate the Bayesian approach we consider the two-dimensional
chapel data set in the left panel of Figure~\ref{fig:Data}, using a
uniform prior for both $\mu=(\cos\varphi,\sin\varphi)$ and $\sigma^2$,
and flat conjugated gamma priors for $\rho_L$ and $\alpha$, see
Figure~\ref{fig:marg_param_hist}.  Our posterior results for $\rho_L$,
$\varphi$, and $\alpha$ were sensitive to the choice of prior
distribution for $\kappa$. For small values of $\kappa$, i.e., values
less than $30$, meaningless posterior results appeared, since
$\varphi$ was approximately uniform, and for $\varphi$ close to zero,
$\rho_L$ tended to zero and hence $\alpha$ tended to infinity. On the
other hand, very large values of $\kappa$ caused a very concentrated
posterior distribution for $\varphi$. As a compromise, after some
experimentation, we fixed $\kappa=40$.

For this model an extension of the observation window $W=[-.5,.5]^2$
to $W_{\text{ext}}=[-0.55,0.55]^2$ seemed large enough to account for
edge effects. For the posterior simulations we used 200,000
iterations, where one iteration consists of updating all the
parameters and the missing data. We considered trace plots, which have
been omitted here, for the parameters and information about the
missing data, indicating that a burn-in of 5000 iterations is
sufficient.
 
There is a clear distinction between the simulated posterior results
for the parameters and the priors, cf.\ the first four panels in
Figure~\ref{fig:marg_param_hist}. The posterior mean of $\varphi$
($115.02^\circ$) is in close agreement with the result of $117^\circ$
found in Section~\ref{s:examplesK}, and $\sigma$ is unlikely to be
larger than $0.02$, indicating that the points are rather close to the
lines and that the choice of $W_{\mathrm{ext}}$ makes sense.  The
final panel in Figure~\ref{fig:marg_param_hist} shows a good agreement
between the number of chapels (110) and the posterior mean of the
intensity $\rho$ (103.7), though there is some uncertainty in the
posterior distribution of $\rho$. Moreover, the posterior means of
$\rho_L$ and $\alpha$ are 12.9 and 8.4, respectively, which combined
with \eqref{e:2} result in the estimate 108.4 for $\rho$.

\begin{figure}
\centering
\includegraphics[width=8cm]{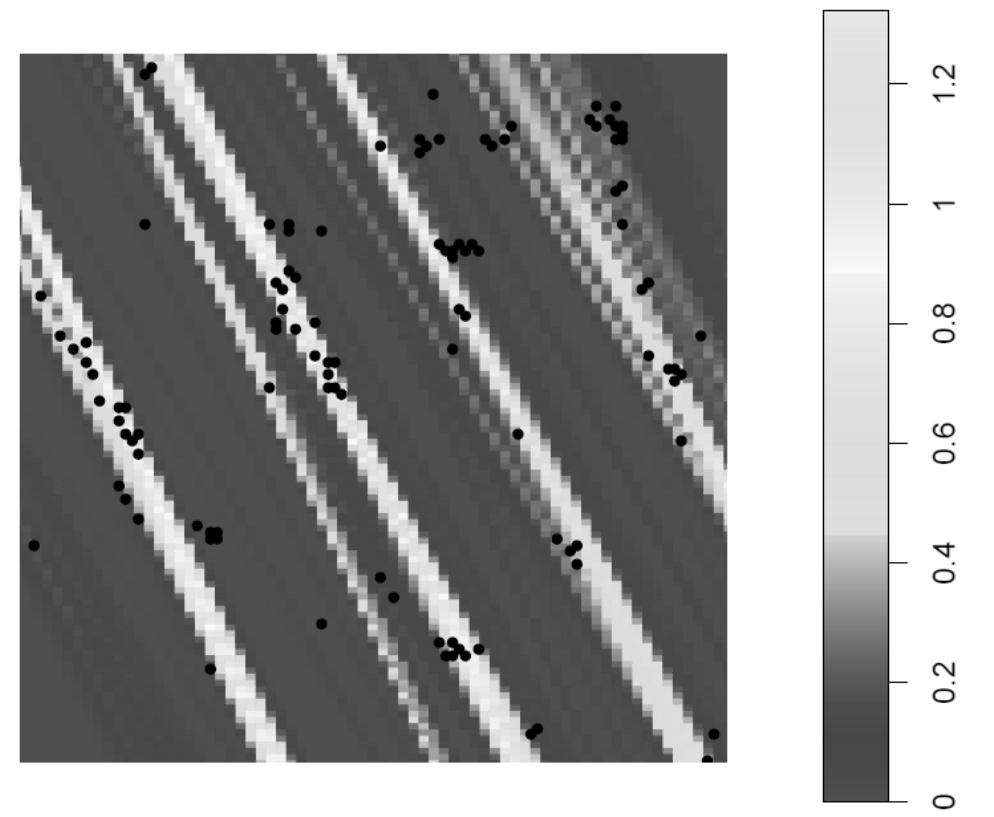}
\caption{Posterior kernel estimate of the density of lines. For comparison, the
  chapel point pattern data set is superimposed.}
\label{fig:lines}
\end{figure}

To illustrate the usefulness of the Bayesian method in detecting
linear structures, Figure~\ref{fig:lines} shows a posterior kernel
estimate of the density of lines within $W$.  The estimate visualizes
where the hidden lines could be, i.e., the lighter areas, and overall
they agree with the point pattern of chapels, which is superimposed in
the figure, though in the upper right corner of the observation window
there is some doubt about whether there should be a single or two
clusters of points.  Specifically, the estimate is obtained from 100
posterior iterations with an equal spacing, and it is the average of
binary pixel representations of the line process, where a pixel has
value 1 if it is intersected by a line, and value 0 otherwise.

\section{Discussion}\label{s:extensions}

\subsection{Choice of structuring element}

Ripley's $K$-function has a ball as structuring element but this is
not useful for detecting anisotropy. Directional $K$-functions have
been suggested using a sector annulus \citep{OhSt-81} or a double cone
\citep{Redenbachetal-09} as the structuring element, while we have
suggested a cylinder.  Another suggestion could be an ellipsoid.

A detailed comparison of the cylindrical $K$-function with the
$K$-function in \citet{Redenbachetal-09} using a double cone as the
structuring element is given in a technical report by F. Safavimanesh
and C. Redenbach from 2016. They conclude that in situations where the
anisotropy is pronounced, a good choice of structuring element may be
important and will depend on the application at hand.  In case of
geometric anisotropy \citep{MoellerHokan-14}, an ellipsoid is
appropriate; when there is a columnar structure as under the Poisson
line cluster point process model or in the data sets considered in
this paper, an elongated cylinder is appropriate; while if regular
point process models are compressed, then a double cone is
appropriate. F. Safavimanesh and C. Redenbach emphasize the importance
of an appropriate choice of the scale and shape parameters used to
specify the structuring element, i.e.\ $r$ and $t$ in case of
$K_{u}(r,t)$. They notice that prior information, e.g., the diameter
of the clusters/mini-columns of points in case of the pyramidal cells
\citep{RSDRMN-15}, can be used to determine interesting ranges of $r$
values.

\subsection{Non-stationary case}

In some applications it is relevant to use a non-constant intensity
function $\rho(x)$. Suppose we assume second order intensity
reweighted stationarity \citep*{baddeley:moeller:waagepetersen:97} and
$g(x)$ still denotes the pair correlation function. This means
intuitively, that if $x_1,x_2\in\mathbb R^d$ are distinct locations and
$B_1,B_2$ are infinitesimally small sets of volumes $\mathrm d
x_1,\,\mathrm dx_2$ and containing $x_1,x_2$, respectively, then
$\rho(x_1)\rho(x_2)g( x_1-x_2)\,\mathrm dx_1\,\mathrm dx_2$ is the
probability for $X$ having a point in each of $B_1$ and $B_2$. Our
definition \eqref{e:Kg} still applies, while \eqref{e:ksum} becomes
\[
K_{u}(r,t)=\frac{1}{|W|} E\left[\sum_{x_1,
    x_2\in X:\,x_1\not=
    x_2}\frac{\mathbb1\{x_1\in W,x_2-x_1\in
    C_{u}(r,t)\}}{\rho(x_1)\rho(
    x_2)}\right],\quad r>0,\ t>0,
\] 
which in turn can be used when deriving non-parametric estimates. 

Assumption (B2) may be relaxed to obtain a non-stationary Poisson line
cluster point process model for $X$, assuming that for each
line $l_i$ the Poisson process $Q_i$ has intensity function
$\alpha_i(x)=\alpha(x)$ for $x\in l_i$, where
$\alpha$ is a non-negative function which is locally integrable on any
line in $\mathbb R^d$.  Then \eqref{e:1} should be replaced by
\begin{equation}\label{e:abcde}
\Lambda(x)=
\sum_{i=1}^\infty
\alpha\{(x\cdot u_i)u_i+y_i\}
k_{u_I^\bot}\{p_{u_i^\bot}(x-y_i)\},
\quad x\in\mathbb R^d.
\end{equation}
However, this non-stationary extension of the model will be harder to
analyze, e.g., moment results as established
in Section~\ref{s:mom} for the stationary case will in general not
easily extend, and it turns out that the model is not second order
intensity reweighted stationary except in a special case discussed
below. Moreover, while our statistical methodology in
Section~\ref{s:mom} can be straightforwardly extended, the Bayesian
computations in Section~\ref{s:Bayesinf} become harder.

In \eqref{e:abcde} assume that the intensity function
$\alpha(y^{(1)},\ldots,y^{(d-1)},x^{(d)})=\alpha(x^{(d)})$ does not
depend on $y=(y^{(d)},\ldots,y^{(d-1)},0)\in H$.  Let
$c=\int_I\alpha\{x^{(d)}\}\,\mathrm dx^{(d)}$.  Then, using a notation as in
Section~\ref{s:mom}, it can be shown that $X$ is second order
intensity reweighted stationary, $X_I$ is still a Neyman--Scott
process, while the $x^{(d)}$-coordinates of the points in $X_W$ are
independent and identically distributed with density $\alpha\{x^{(d)}\}/c$,
and they are independent of $X_W$. Therefore statistical inference
simply splits into modelling the density $\alpha\{x^{(d)}\}/c$ based on the
$x^{(d)}$-coordinates of the points in $X_W$ and inferring
$(c,\rho_L,\sigma^2)$ by considering $X_I$ along similar lines as in
Section~\ref{s:mom} but with $\alpha\{x^{(d)}\}|I|$ replaced by $c$.

\subsection*{Acknowledgments}

Supported by the Danish Council for Independent Research | Natural Sciences,
and by the Centre for Stochastic Geometry and Advanced Bioimaging,
funded by
the Villum Foundation. We thank Jens Randel Nyengaard, Karl-Anton
Dorph-Petersen, and Ali H.\ Rafati for collecting the
three-dimensional pyramidal cell data set.
\section*{Appendix A: Intensity and rose of direction for a Poisson
  line process}

First we give the definition of the intensity $\rho_L$ and the rose of
direction $\mathcal R$ for a general stationary line process
$L=\{l_1,l_2,\ldots\}$ in $\mathbb R^d$.  Let $|\cdot|_1$ denote
one-dimensional Lebesgue measure, $\mathrm dt$ Lebesgue measure on the
real line, $A\subseteq\mathbb R^d$ an arbitrary Borel set with volume
$|A|\in(0,\infty)$, and an $B\subseteq\mathbb S^{d-1}$ arbitrary Borel
set. Then by definition and since $L$ is stationarity,
\begin{equation}
\rho_L = E\left(\sum_{i=1}^\infty
|l_i\cap A|_1/|A|\right)
\label{aaaa}
\end{equation}
does not depend on the choice of $A$, and provided $0<\rho_L<\infty$,
\begin{equation}
\mathcal R(B)= E\left\{\sum_{i=1}^\infty
|l_i\cap A|_1\mathrm1(u_i\in B)/\left(\rho_L|A|\right)\right\}\label{a}
\end{equation}
does not depend on the choice of $A$ and is seen to be a probability
measure.

Second we assume that $L$ is a stationary Poisson line
process as in Section~\ref{s:defprop}. Then
\begin{align}
E\left\{\sum_{i=1}^\infty
|l_i\cap A|_1\mathrm1(u\in B)\right\}=
\,& E\left\{\sum_{i=1}^\infty
\int\mathrm1(y_i+tu_i\in
A,\, u_i\in B)\,\mathrm dt\right\}\label{b}\\\
=\,&\beta\int\int\int\mathrm1(y+tu\in
A,\, u\in B)\,\mathrm dt\,
\lambda(\mathrm dy)\,M(\mathrm du)\label{c}\\
=\,&\beta|A|\int_B1/|u^{(d)}|\,M(\mathrm du)\label{d}.
\end{align}
Here 
\eqref{b} follows from the phase representation of $L$ (see
Section~\ref{s:defprop}), 
\eqref{c} from  the Slivnyak--Mecke theorem for the 
Poisson process $\Phi$ (see e.g.\ \citet{MW-04}), 
and \eqref{d} since $|u^{(d)}|$ is the Jacobian
of the mapping $(t,y)\mapsto y+tu$ with
$(t,y)\in\mathbb R\times H$.
When $B=\mathbb S^{d-1}$ we obtain from \eqref{aaaa} and \eqref{d} 
the first equation in \eqref{e:L}. This together with $0<\beta<\infty$
and $0<\int_{\mathbb S^{d-1}}1/|u^{(d)}|\,M(\mathrm du)<\infty$
imply that $0<\rho_L<\infty$. Thereby the second first equation in
\eqref{e:L} follows for any Borel set $B\subseteq\mathbb S^{d-1}$.

\section*{Appendix B: Moment results for a Poisson line cluster point process}


For Section~\ref{s:mom} it remains to verify \eqref{e:2}--\eqref{e:4}.

{\it Proof of \eqref{e:2}:} By \eqref{e:1}, \eqref{e:1a}, and the
Slivnyak--Mecke theorem for the Poisson process,
\begin{equation}\label{e:3}
  \rho=\alpha\beta\int\int k_{u}\{-p_{u}(y)\}
  \,\lambda(\mathrm dy)\,M(\mathrm du).
\end{equation}
Let $I_d$ be the $d\times d$ identity matrix, and $o_{d-1}$
the origin in $\mathbb R^{d-1}$. For $u\in\mathbb{S}^{d-1}$,
let $A(u)=vv^T$ where $v$ is the
subvector consisting of the first $d-1$ coordinates of $u$ and
$^T$ denotes transpose of a vector or matrix. The Jacobian of the
linear transformation
\[
p_{u}(y)=\left(I_d-uu^T\right)y,\quad
y\in H,
\]
is the square root of the determinant of the $(d-1)\times(d-1)$ matrix 
\begin{align*}
Q&=\left\{\left(I_d-uu^T\right) 
\left(I_{d-1}\ o_{d-1}\right)^T\right\}^T
\left(I_d-uu^T\right) 
\left(I_{d-1}\ o_{d-1}\right)^T\\
&=\left(I_{d-1}\ o_{d-1}\right)
\left(I_d-uu^T\right)
\left(I_{d-1}\ o_{d-1}\right)^T
\\ &=I_{d-1}-A(u).
\end{align*}
Since $A(u)$ is symmetric of rank at most one and has trace
$\mbox{tr}\{A(u)\}=\{u^{(1)}\}^2+\ldots+\{u^{(d-1)}\}^2=1-\{u^{(d)}\}^2$, the
determinant of $Q$ is $1-\mbox{tr}\{A(u)\}=\{u^{(d)}\}^2$.  Combining
this with \eqref{e:3} we obtain
\[\rho=\alpha\beta\int\int k_{u}(-y)/|u^{(d)}|
\,\lambda_{u}(\mathrm dy)\,M(\mathrm du). \]
Thereby 
\eqref{e:2} easily follows from the first identity in \eqref{e:L}. 

{\it Proof of \eqref{e:4}:} By \eqref{e:1} and \eqref{e:1a},
\begin{align}
\rho^2g(x)&=\alpha^2
E\left[\sum_{i\not= j}k_{u_i^\bot}\{p_{u_i^\bot}(-y_i)\}
k_{u_j^\bot}\{p_{u_j^\bot}(x-y_j)\}\right]\nonumber\\
&+
\alpha^2 E\left[\sum_ik_{u_i^\bot}\{p_{u_i^\bot}(-y_i)\}
k_{u_i^\bot}\{p_{u_i^\bot}(x-y_i)\}\right]\nonumber\\
&=\rho^2+\alpha^2\beta\int k_{u^\bot}\{p_{u^\bot}(-y)\}
k_{u^\bot}\{p_{u^\bot}(x-y)\}
\,\lambda(\mathrm dy)\,M(\mathrm du)\label{e:5}
\end{align}
using the extended
Slivnyak--Mecke theorem for the Poisson process and
the proof of \eqref{e:2} to obtain that the first expectation is equal
to $\rho^2$, and the Slivnyak--Mecke theorem for the Poisson process 
to obtain that the second expectation is equal to the last term.
Combining \eqref{e:2} and \eqref{e:5} with the result 
for the Jacobian considered above,
we obtain \eqref{e:4}. 

\section*{Appendix C: Hybrid Markov chain Monte Carlo algorithm}
This appendix details the Markov chain Monte Carlo algorithm for the
Bayesian approach considered in Section~\ref{s:Bayesinf}, with
independent prior densities for the parameters and posterior density
given by \eqref{e:post}--\eqref{(20)}.  As in
Section~\ref{s:Bayesinf}, we consider conjugated gamma densities
$p(\alpha)$ and $p(\rho_L)$, and denote their shape parameters by
$a_1$ and $a_2$ and their inverse scale parameters by $b_1$ and $b_2$,
respectively. The remaining parameters have no (well-known) conjugate
priors, cf.\ \eqref{e:den1} and \eqref{e:den2}, and thus we consider
generic prior densities $p(\mu)$, $p(\kappa)$, and $p(\sigma^2)$.

In each iteration of the Markov chain Monte Carlo algorithm we update
first each of the parameters and second the missing data.  We use a
Gibbs update for $\alpha$ respective $\rho_L$, noting that the
conditional distribution of $\alpha$ given the rest is a gamma
distribution with shape parameter $a_1+n$ and inverse scale parameter
$b_1+\int_W\sum_{j=1}^k f\{p_{u_j^\bot}(x-y_j)\mid\sigma^2\}\, \mathrm
dx$, and the conditional distribution of $\rho_L$ given the rest is a
gamma distribution with shape parameter $a_2+k$ and inverse scale
parameter $b_2+I(\mu,\kappa)$.  Below we describe the individual
proposals and Hastings ratios for the remaining parameters and the
missing data (in the case of Section~\ref{s:Bayesinf} where the value
of $\kappa$ is fixed, we can of course just ignore the update of
$\kappa$ described below).  As usual, for each type of update, the
proposal is accepted with probability $\min\{1,R\}$, where $R$ is the
corresponding Hastings ratio.  We denote
$(\rho_L,\mu,\kappa,\alpha,\sigma^2,\{(y_1,u_1),\ldots, (y_k,u_k)\})$
the current state of the algorithm, where $n\ge 1$, $k\ge1$ (since
$k=0$ implies $n=0$, which is not a case of interest), and
$l(y_i,u_i)\cap W_{\text{ext}}\not=\emptyset, i=1,\ldots,k$.

For 
$\mu$, $\kappa$, and $\sigma^2$,
we use Metropolis random walk updates, with 
a von Mises--Fisher 
proposal $\mu'\sim f(\cdot\mid\mu,\kappa_0)$
and
normal proposals $\kappa'\sim N(\kappa,\sigma_{0,1}^2)$ 
and $\sigma'^2\sim N(\sigma^2,\sigma_{0,2}^2)$, where
$\kappa_0,\sigma_{0,1}^2,\sigma_{0,2}^2>0$ are  
tuned so that the mean acceptance probabilities are
between 20--45\%  (as recommended in \citet*{Robertsetal-97}). 
The Hastings
ratios for the acceptance probabilities are
\[R_{\mu}=\frac{p(\mu')}{p(\mu)}
\exp\left[\rho_L\left\{I(\mu,\kappa)-
I(\mu',\kappa)\right\}\right]\prod_{j=1}^k\frac
{f(u_j\mid\mu',\kappa)}
{f(u_j\mid\mu,\kappa)},\]
\[R_\kappa=\mathbb1(\kappa'>0)\frac{p(\kappa')}{p(\kappa)}
\exp\left[\rho_L\left\{I(\mu,\kappa)-
I(\mu,\kappa')\right\}\right]\prod_{j=1}^k\frac
{f(u_j\mid\mu,\kappa')}
{f(u_j\mid\mu,\kappa)},\]
and
\begin{align*}
  R_{\sigma^2}=&\,\mathbb1(\sigma'^2>0)\frac{p(\sigma'^2)}{p(\sigma^2)}\exp\Bigg(\alpha
  \sum_{j=1}^k\left[\int_W
    f\{p_{u_j^\bot}(x-y_j)\mid\sigma^2\}\,\mathrm dx \right.\\ &
  \left.  -\int_W f\{p_{u_j^\bot}(x-y_j)\mid\sigma'^2\}\,\mathrm
    dx\right]\Bigg) \prod_{i=1}^n
  \frac{\sum_{j=1}^kf\{p_{u_j^\bot}(x_i-y_j)\mid\sigma'^2\}}
  {\sum_{j=1}^kf\{p_{u_j^\bot}(x_i-y_j)\mid\sigma^2\}}.
\end{align*}
For $R_{\sigma^2}$ each integral is calculated by a simple Monte Carlo method after
making a change of variables from $x$ to its Cartesian coordinates
in a system centered at $y_j$ and with axes given by
$u_j$ and $u_j^{\bot}$.

For the missing data, we adapt the birth-death-move Metropolis-Hastings 
algorithm in \citet{GeyerMoeller-94} as follows. 
Each of the birth/death/move proposals happens with
 probability  $1/3$ and consists of the following action.      
A birth proposal is the proposal of adding a new point
$(y,u)$, where 
$u\sim f(\cdot\mid\mu,\kappa)$ and 
$y$ conditional on 
$u$ is uniformly distributed on $J_{u}$. Then, as
explained below, the
Hastings ratio is
\begin{align}\label{e:Rbirth}
R_{\text{birth}}=&\,\frac{\rho_L\lambda(J_{u})|u^{(d)}|}{k+1}
\mathbb1\{l(y,u)\cap W_{\text{ext}}\not=\emptyset\}\nonumber\\
&\times\exp\left[-\alpha\int_W f\{p_{u^\bot}(x-y)\mid\sigma^2\}
\mathrm dx\right]\prod_{i=1}^n
\left[1+\frac{f\{p_{u^\bot}(x_i-y)\mid\sigma^2\}}
{\sum_{j=1}^kf\{p_{u_j^\bot}(x_i-y_j)\mid\sigma^2\}}\right].
\end{align}
To stress the dependence on 
 $\{(y_1,u_1),\ldots,(y_k,u_k)\}$  
and $(y,u)$, 
write 
$R_{\text{birth}} = R_{\text{birth}}(y_1,u_1,\ldots,y_k,u_k;y,u)$
  (obviously, it also depends on 
$\rho_L,\alpha,\sigma^2$, and $\{x_1,\ldots,x_n\}$). 
A death proposal is the proposal of generating a uniform
$j\in\{1,\ldots,k\}$ (provided $k>1$; if $k=1$, we do nothing and keep
the current state) and deleting 
$(y_j,u_j)$. Then, as
explained below, 
the Hastings ratio is
\begin{equation}\label{e:Rdeath}
R_{\text{death}}=1/R_{\text{birth}}(y_1,u_1,
\ldots,y_{j-1},u_{j-1},y_{j+1},
u_{j+1}\ldots,y_k,u_k;y_j,u_j).
\end{equation}
Finally, 
a move proposal is the proposal of selecting a uniform $j\in\{1,\ldots,k\}$
and replacing $(y_j,u_j)$ by 
$(y'_j,u'_j)$, where 
$u'_j\sim f(\cdot\mid\mu,\kappa)$ and 
$y'_j$ conditional on
$u'_j$ is uniformly distributed on $J_{u'_j}$. Since
this can be considered as first a death proposal and second a birth
proposal, the
Hastings ratio is
\begin{equation*}\label{e:Rmove}
R_{\text{move}}=\frac{R_{\text{birth}}(y_1,u_1,
\ldots,y_{j-1},u_{j-1},y_{j+1},
u_{j+1}\ldots,y_k,u_k;y'_j,u'_j)}
{R_{\text{birth}}(y_1,u_1,\ldots,y_{j-1},
u_{j-1},y_{j+1},u_{j+1}\ldots,y_k,u_k;y_j,u_j)}.
\end{equation*}

It remains to explain how we obtained the Hastings
ratios~(\ref{e:Rbirth})--(\ref{e:Rdeath}), where we notice the
following facts. 
The reference Poisson process $\Phi_{0,S}$ has intensity
 measure 
\[\zeta(\mathrm dy,\mathrm du)=|u^{(d)}|\frac{\Gamma(d/2)}{2\pi^{d/2}}\lambda(\mathrm dy)\nu_{d-1}(\mathrm du). \] 
Further, conditional on the data $X_W=\{x_1,\ldots,x_n\}$ and the parameters $\rho_L$, $\mu$, $\kappa$, $\alpha$, $\sigma^2$, the target process $\Phi_{S}$ has density 
\begin{align*}
&f[\{(y_1,u_1),\ldots,
(y_k,u_k)\}\mid\{x_1,\ldots,x_n\},\alpha,{\sigma}^2,\rho_L,\mu,\kappa]\\
&\propto f[\{x_1,\ldots,x_n\}\mid
\{(y_1,u_1),\ldots,
(y_k,u_k)\},\alpha,\sigma^2]
f[\{(y_1,u_1),\ldots,
(y_k,u_k)\}\mid\rho_L,\mu,\kappa]\label{e:likelihood}
\end{align*}
with respect to the distribution of $\Phi_{0,S}$. Furthermore, if 
a birth $(y,u)$ is proposed, then it has density 
\[
f(y,u\mid\mu,\kappa)=\frac{f(u\mid\mu,\kappa)\mathbb1(y\in
  J_{u})/\lambda(J_{u})}{|u^{(d)}|\Gamma(d/2)/(2\pi^{d/2})}
\] 
with respect to $\zeta$.
Consequently, by \citet{GeyerMoeller-94}, the Hastings ratio for the proposed birth is 
\begin{align*}
R_{\text{birth}}=&\,\frac{f[\{(y_1,u_1),\ldots,
(y_k,u_k),(y,u)\}\mid\{x_1,\ldots,x_n\},\alpha,{\sigma}^2,\rho_L,\mu,\kappa]}{f[\{(y_1,u_1),\ldots,
(y_k,u_k)\}\mid\{x_1,\ldots,x_n\},\alpha,{\sigma}^2,\rho_L,\mu,\kappa]}\times\frac{1/(k+1)}{f(y,u\mid\mu,\kappa)}\\
=&\,\frac{2\pi^{d/2}}{\Gamma(d/2)}\rho_L f(u\mid\mu,\kappa)\mathbb1(y\in J_{u})\frac{1}{(k+1)f(y, u\mid\mu,\kappa)}\\
&\times\exp\left[-\alpha\int_W f\{p_{u^\bot}(x-y)\mid\sigma^2\}\mathrm dx\right]\prod_{i=1}^n
\left[1+\frac{f\{p_{u^\bot}(x_i-y)\mid\sigma^2\}}
{\sum_{j=1}^kf\{p_{u_j^\bot}(x_i-y_j)\mid\sigma^2\}}\right]
\end{align*}
which is equal to \eqref{e:Rbirth}. Thereby, refering again to \citet{GeyerMoeller-94}, we obtain \eqref{e:Rdeath}. 

\bibliographystyle{biometrika}
\bibliography{cylK} 

\begin{thebibliography}{29}
\expandafter\ifx\csname natexlab\endcsname\relax\def\natexlab#1{#1}\fi

\bibitem[{Baddeley et~al.(2000)Baddeley, M{\o}ller \&
  Waagepetersen}]{baddeley:moeller:waagepetersen:97}
\textsc{Baddeley, A.}, \textsc{M{\o}ller, J.} \& \textsc{Waagepetersen, R.}
  (2000).
\newblock Non- and semi-parametric estimation of interaction in inhomogeneous
  point patterns.
\newblock \textit{Statistica Neerlandica} \textbf{54}, 329--350.

\bibitem[{Baddeley \& Turner(2005)}]{spatstat}
\textsc{Baddeley, A.} \& \textsc{Turner, R.} (2005).
\newblock {Spatstat}: an {{\sf R}} package for analyzing spatial point
  patterns.
\newblock \textit{Journal of Statistical Software} \textbf{12}, 1--42.

\bibitem[{Chiu et~al.(2013)Chiu, Stoyan, Kendall \& Mecke}]{Stoyanetal-94}
\textsc{Chiu, S.~N.}, \textsc{Stoyan, D.}, \textsc{Kendall, W.~S.} \&
  \textsc{Mecke, J.} (2013).
\newblock \textit{Stochastic Geometry and Its Applications}.
\newblock John Wiley and Sons, Chichester, 3rd ed.

\bibitem[{Diggle et~al.(1995)Diggle, Chetwynd, H{\"a}ggkvist \&
  Morris}]{Diggleetal-95}
\textsc{Diggle, P.~J.}, \textsc{Chetwynd, A.}, \textsc{H{\"a}ggkvist, R.} \&
  \textsc{Morris, S.} (1995).
\newblock Second order analysis of space-time clustering.
\newblock \textit{Statistical Methods in Medical Research} \textbf{4},
  124--136.

\bibitem[{Diggle \& Gratton(1984)}]{mincontrast-84}
\textsc{Diggle, P.~J.} \& \textsc{Gratton, R.} (1984).
\newblock Monte {C}arlo methods of inference for implicit statistical models
  (with discussion).
\newblock \textit{Journal of the Royal Statistical Society: Series B
  (Statistical Methodology)} \textbf{46}, 193--212.

\bibitem[{Gabriel \& Diggle(2009)}]{GabrielDiggel-09}
\textsc{Gabriel, E.} \& \textsc{Diggle, P.~J.} (2009).
\newblock Second-order analysis of inhomogeneous spatio-temporal point process
  data.
\newblock \textit{Statistica Neerlandica} \textbf{63}, 43--51.

\bibitem[{Geyer \& M{\o}ller(1994)}]{GeyerMoeller-94}
\textsc{Geyer, C.} \& \textsc{M{\o}ller, J.} (1994).
\newblock Simulation procedures and likelihood inference for spatial point
  processes.
\newblock \textit{Scandinavian Journal of Statistics} \textbf{21}, 359--373.

\bibitem[{Gilks et~al.(1996)Gilks, Richardson \& Spiegelhalter}]{Gilksetal-96}
\textsc{Gilks, W.~R.}, \textsc{Richardson, S.} \& \textsc{Spiegelhalter, D.~J.}
  (1996).
\newblock \textit{Markov Chain Monte Carlo in Practice}.
\newblock Chapman and Hall, London.

\bibitem[{Guan(2006)}]{Guan-06}
\textsc{Guan, Y.} (2006).
\newblock A composite likelihood approach in fitting spatial point process
  models.
\newblock \textit{Journal of the American Statistical Association}
  \textbf{101}, 1502--1512.

\bibitem[{Guan et~al.(2006)Guan, Sherman \& Calvin}]{guan-etal:06}
\textsc{Guan, Y.}, \textsc{Sherman, M.} \& \textsc{Calvin, J.~A.} (2006).
\newblock Assessing isotropy for spatial point processes.
\newblock \textit{Biometrics} \textbf{62}, 119--125.

\bibitem[{Illian et~al.(2008)Illian, Penttinen, Stoyan \&
  Stoyan}]{Illianetal-08}
\textsc{Illian, J.}, \textsc{Penttinen, A.}, \textsc{Stoyan, H.} \&
  \textsc{Stoyan, D.} (2008).
\newblock \textit{Statistical Analysis and Modelling of Spatial Point
  Patterns}.
\newblock John Wiley and Sons, New York.

\bibitem[{M{\o}ller \& Toftaker(2014)}]{MoellerHokan-14}
\textsc{M{\o}ller, J.} \& \textsc{Toftaker, H.} (2014).
\newblock Geometric anisotropic spatial point pattern analysis and {C}ox
  processes.
\newblock \textit{Scandinavian Journal of Statistics} \textbf{41}, 414--435.

\bibitem[{M{\o}ller \& Waagepetersen(2004)}]{MW-04}
\textsc{M{\o}ller, J.} \& \textsc{Waagepetersen, R.} (2004).
\newblock \textit{Statistical Inference and Simulation for Spatial Point
  Processes}.
\newblock Chapman and Hall/CRC, Boca Raton.

\bibitem[{M{\o}ller \& Waagepetersen(2007)}]{MW-07}
\textsc{M{\o}ller, J.} \& \textsc{Waagepetersen, R.} (2007).
\newblock Modern statistics for spatial point processes (with discussion).
\newblock \textit{Scandinavian journal of Statistics} \textbf{34}, 643--711.

\bibitem[{M{\o}ller \& Waagepetersen(2016)}]{moeller:waagepetersen:16}
\textsc{M{\o}ller, J.} \& \textsc{Waagepetersen, R.} (2016).
\newblock Some recent developments in statistics for spatial point pattern
  analysis.
\newblock \textit{Annual Review of Statistics and Its Application} (submitted
  invited paper).

\bibitem[{Mountcastle(1957)}]{Mountcastle-57}
\textsc{Mountcastle, V.~B.} (1957).
\newblock Modality and topographic properties of single neurons of cat’s
  somatic sensory cortex.
\newblock \textit{Journal of Neurophysiology} \textbf{20}, 408--434.

\bibitem[{Mugglestone \& Renshaw(1996)}]{MuggRen-96}
\textsc{Mugglestone, M.} \& \textsc{Renshaw, E.} (1996).
\newblock A practical guide to the spectral analysis of spatial point
  processes.
\newblock \textit{Computational Statistics {$\&$} Data Analysis} \textbf{21},
  43--65.

\bibitem[{Myllym{\"a}ki et~al.(2016)Myllym{\"a}ki, Mrkvi{\v c}ka, Seijo \&
  Grabarnik}]{Myllymaki-16}
\textsc{Myllym{\"a}ki, M.}, \textsc{Mrkvi{\v c}ka, T.}, \textsc{Seijo, H.} \&
  \textsc{Grabarnik, P.} (2016).
\newblock Global envelope tests for spatial processes.
\newblock \textit{Journal of the Royal Statistical Society: Series B
  (Statistical Methodology)} (to appear).

\bibitem[{Nicolis et~al.(2010)Nicolis, Mateu \& D'Ercole}]{Nicolisetal-10}
\textsc{Nicolis, O.}, \textsc{Mateu, J.} \& \textsc{D'Ercole, R.} (2010).
\newblock Testing for anisotropy in spatial point processes.
\newblock In \textit{Proceedings of the Fifth International Workshop on
  Spatio-Temporal Modelling (METMA5)}, G.-M. et~al., ed. Unidixital, Santiago
  de Compostela.

\bibitem[{Ohser \& Stoyan(1981)}]{OhSt-81}
\textsc{Ohser, J.} \& \textsc{Stoyan, D.} (1981).
\newblock On the second-order and orientation analysis of planar stationary
  point processes.
\newblock \textit{Biometrical Journal} \textbf{23}, 523–533.

\bibitem[{Rafati et~al.(2016)Rafati, Safavimanesh, Dorph-Petersen, Rasmussen,
  M{\o}ller \& Nyengaard}]{RSDRMN-15}
\textsc{Rafati, A.}, \textsc{Safavimanesh, F.}, \textsc{Dorph-Petersen, K.},
  \textsc{Rasmussen, J.~G.}, \textsc{M{\o}ller, J.} \& \textsc{Nyengaard,
  J.~R.} (2016).
\newblock Detection and spatial characterization of minicolumnarity in the
  human cerebral cortex.
\newblock \textit{Journal of Microscopy} \textbf{261}, 115--126.

\bibitem[{Redenbach et~al.(2009)Redenbach, S{\"a}rkk{\"a}, Freitag \&
  Schladitz}]{Redenbachetal-09}
\textsc{Redenbach, C.}, \textsc{S{\"a}rkk{\"a}, A.}, \textsc{Freitag, J.} \&
  \textsc{Schladitz, K.} (2009).
\newblock Anisotropy analysis of pressed point processes.
\newblock \textit{Advances in Statistical Analysis} \textbf{93}, 237–261.

\bibitem[{Ripley(1976)}]{ripley:76}
\textsc{Ripley, B.~D.} (1976).
\newblock The second-order analysis of stationary point processes.
\newblock \textit{Journal of Applied Probability} \textbf{13}, 255--266.

\bibitem[{Ripley(1977)}]{ripley:77}
\textsc{Ripley, B.~D.} (1977).
\newblock Modelling spatial patterns (with discussion).
\newblock \textit{Journal of the Royal Statistical Society: Series B
  (Statistical Methodology)} \textbf{39}, 172--212.

\bibitem[{Roberts et~al.(1997)Roberts, Gelman \& Gilks}]{Robertsetal-97}
\textsc{Roberts, G.~O.}, \textsc{Gelman, A.} \& \textsc{Gilks, W.~R.} (1997).
\newblock Weak convergence and optimal scaling of random walk {M}etropolis
  algorithms.
\newblock \textit{Annual of Applied Probability} \textbf{7}, 110–120.

\bibitem[{Rosenberg(2004)}]{Rosenberg-04}
\textsc{Rosenberg, M.~S.} (2004).
\newblock Wavelet analysis for detecting anisotropy in point patterns.
\newblock \textit{Journal of Vegetation Science} \textbf{15}, 277–284.

\bibitem[{Stoyan(1991)}]{Stoyan-91}
\textsc{Stoyan, D.} (1991).
\newblock Describing the anisotropy of marked planer point process.
\newblock \textit{Statistics: A Journal of Theoretical and Applied Statistics}
  \textbf{22}, 449--462.

\bibitem[{Stoyan \& Bene{\v s}(1991)}]{StoyanBenes-91}
\textsc{Stoyan, D.} \& \textsc{Bene{\v s}, V.} (1991).
\newblock Anisotropy analysis for particle systems.
\newblock \textit{Journal of Microscopy} \textbf{164}, 159--168.

\bibitem[{Stoyan \& Stoyan(1995)}]{StoyStoy-94}
\textsc{Stoyan, D.} \& \textsc{Stoyan, H.} (1995).
\newblock \textit{Fractals, Random Shapes and Point Fields}.
\newblock John Wiley and Sons, Chichester.

\end{thebibliography}

\end{document}